\documentclass[12pt]{article}

\usepackage{amssymb,amsmath}

\usepackage{pdfpages}

\newtheorem{theorem}{Theorem}[section]
\newtheorem{corollary}[theorem]{Corollary}

\newtheorem{proposition}[theorem]{Proposition}

\title{On Shilnikov's scenario in 3D: Topological chaos
for vectorfields of class $C^1$}

\author{Hans-Otto Walther}

\begin{document}
	
\maketitle

\section{Introduction}

In his seminal paper \cite{S3} Shilnikov considered the differential equation  
\begin{equation}
x'(t)=V(x(t))\in\mathbb{R}^3
\end{equation}
for $V$ analytic with $V(0)=0$ so that there is a homoclinic solution $h_V:
\mathbb{R}\to\mathbb{R}^3$, $0\neq h_V(t)\to0$ for $|t|\to\infty$,
and the eigenvalues $u\in\mathbb{R}$ and $\sigma\pm i\,\mu\in\mathbb{C}$ of the derivative $DV(0)$ satisfy $u>0$, $\sigma<0<\mu$, with

\medskip

(H)$\quad\quad 0<\sigma+u.$

\medskip

Shilnikov's 
result in \cite{S3} is that close to the homoclinic orbit $h_V(\mathbb{R})$ there exist countably many periodic orbits.

\medskip

A related, slightly stronger statement about complicated motion is conjugacy with the shift $(s_n)_{n=-\infty}^{\infty}\mapsto(s_{n+1})_{n=-\infty}^{\infty}$
in two symbols $s_n\in\{0,1\}$, for a return map, 
which is given by intersections of solutions with a transversal to the homoclinic orbit. See the monographs \cite{GH,Wi1,Wi2,HSD} for work on the verification of this property and of related ones. A detailed presentation for $V$ linear close to equilibrium is contained in \cite{HSD}.

\medskip

Shilnikov-type results on complicated motion close to a homoclinic orbit have also been obtained for semiflows in infinite-dimensional spaces, e. g. in \cite{W1,LWW}
on delay differential equations which are linear close to the stationary state. These results are related to Shilnikov's work on vectorfields on $\mathbb{R}^4$ \cite{S4} which is analogous to \cite{S3} except for the presence of pairs of complex conjugate eigenvalues of $DV(0)$ in either halfplane.

\medskip

The precursor \cite{W2} of the present paper contains a proof that for $V$ in Eq. (1) twice continuously differentiable there exists {\it topological chaos} close to the homoclinic orbit, in the sense that for any given symbol sequence $(s_n)_{n=-\infty}^{\infty}$
there is a trajectory $(x_n)_{n=-\infty}^{\infty}$ of the return map taking values in disjoint sets $M_0,M_1$ according to $x_n\in M_{s_n}$.

\medskip

In the present paper we prove existence of topological chaos
for vectorfields $V$ which only need to be once continuously differentiable - a kind of minimal smoothness for Shilnikov's scenario. See Proposition 8.1 and
Theorem 8.2 for precise statements of the main results.

\medskip

In contrast to the approach in \cite{W2} we work solely with flows and flowlines, not with solutions of underlying differential equations. The reason is simply that transforming the flow $F_V$
of Eq. (1) by means of a diffeomorphism preserves smoothness whereas transforming Eq. (1) as in \cite{W2}, into a better tractable differential equation with flat local stable and unstable  manifolds, would yield a vectorfield which is in general only continuous and thereby not good enough for arguments as used in \cite{W2}. The appendix in Section 9 below describes how to transform  $F_V$ for $V:\mathbb{R}^3\to\mathbb{R}^3$ continuously differentiable and, say, bounded,
into  a continuously differentiable flow $F:\mathbb{R}\times\mathbb{R}^3\to\mathbb{R}^3$ with a homoclinic flowline and flat local stable and unstable manifolds, and with further properties, altogether listed as (F1)-(F5) in Section 2 below. For convenience we turn to scaled flows given by $F_{\epsilon}(t,x)=\frac{1}{\epsilon}F(t,\epsilon\,x), \,\,\epsilon>0$, all of which are equivalent to $F$. In the sequel we investigate the behaviour of these scaled flows inside and outside of a fixed neighbourhood of the origin, instead of studying $F$ with respect to a family of  shrinking neighbourhoods. The content of Sections 2-6 parallels  its
counterpart in \cite{W2}, but working with the flows instead of a differential equation necessitates modifications, among them in Section 3 another access to angles along projections of flowlines into the stable plane, and a rearrangement of arguments in Sections 4-6.

\medskip

As in \cite{W2} (and following Shilnikov \cite{S3}) we introduce a return map, now only for particular small $\epsilon>0$, with domain in a transversal to the homoclinic flowline, on one side of the flat local stable manifold, as shown at the top of Figure 1. 
\begin{figure}
	\includegraphics[page=1,scale=0.7]{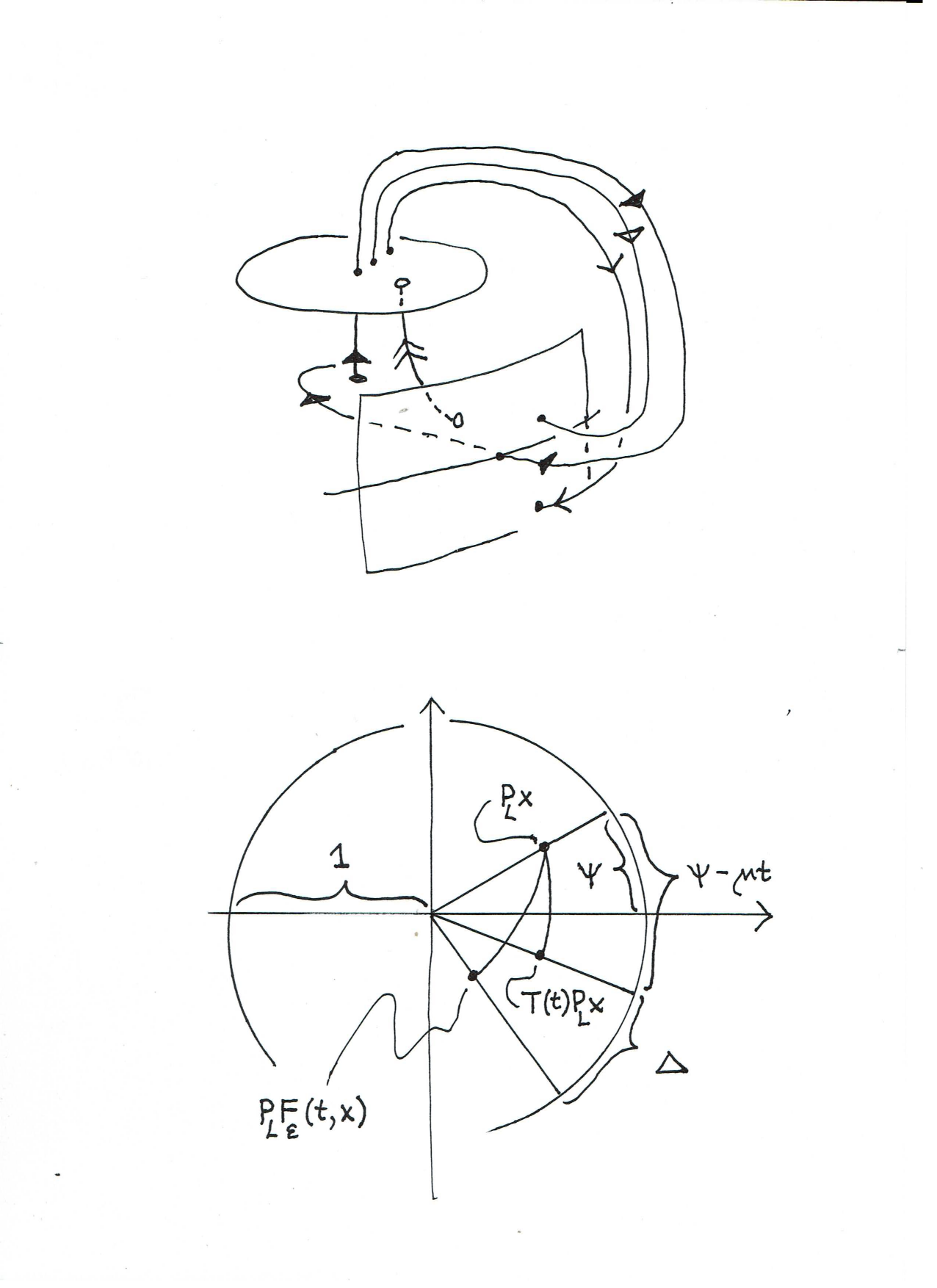}
	\caption{Top: The return map as a composition of the inner map with the exterior map. Bottom: Angles from Proposition 3.1}
\end{figure}
Expressed in suitable coordinates the return map becomes a map from a rectangle into the plane. Section 7 shows how this map turns curves which connect certain horizontal levels in the rectangle into spirals around the origin. This suffices for the proof of Proposition 8.1 about one-directional topological chaos (with forward symbol sequences $(s_n)_{n=0}^{\infty}$), and avoids the verification of covering relations for 2-dimensional subsets of the rectangle. Theorem 8.2 extends the result of Proposition 8.1 to entire
symbol sequences $(s_n)_{n=-\infty}^{\infty}$, by means of familiar compactness
arguments. Sections 7 and 8 proceed almost exactly as in \cite{W2}. We include all details of proofs also in these sections in order to keep  the paper be self-contained.

\medskip

The choice of $\Delta_2$ in Section 7 shows that actually we obtain a countable family of sets of complicated trajectories of the return map. What remains open, among others, is existence of periodic orbits corresponding to periodic symbol sequences. Also of interest might be a version of the present approach for Shilnikov's scenario
in $\mathbb{R}^4$ \cite{S4}.

\medskip

{\bf Notation, preliminaries}. A forward trajectory of a map $f:M\supset\,dom\to M$ is a sequence $(x_j)_0^{\infty}$ in $dom$ with $x_{j+1}=f(x_j)$ for all integers $j\ge0$. 
Entire trajectories are defined analogously, with all integers as indices.

For a vectorspace $X$, $x\in X$, and $M\subset X$, we set $x\pm M=\{y\in X:y\mp x\in M\}$. Similarly, for $A\subset\mathbb{R}$ and $x\in X$, $Ax=\{y\in X:\mbox{For some}\,\, a\in A, y=ax\}$. 

The interior, the boundary, and the closure of a subset of a topological space are denoted by $int\,M, \partial \,M,$ and $cl\,M,$ respectively.

A curve is a continuous map from an interval $I\subset\mathbb{R}$ into a topological space.

Components of vectors in Euclidean spaces $\mathbb{R}^n$ are indicated by lower indices.The inner product on $\mathbb{R}^n$ is written as $<x,y>=\sum_{i=1}^nx_iy_i$, and we use the Euclidean norm given by $|x|=\sqrt{<x,x>}$. The vectors of the canonical orthonormal basis on $\mathbb{R}^n$ are denoted by $e_j$, $j=1,\ldots,n$, $e_{j,j}=1$ and $e_{j,k}=0$ for $j\neq k$. In $\mathbb{R}^3$ we write $L=\mathbb{R}e_1\oplus\mathbb{R}e_2$ and $U=\mathbb{R}e_3$. The  associated projections $\mathbb{R}^3\to\mathbb{R}^3$ onto $L$ and onto $U$ are denoted by $P_L$ and $P_U$, respectively. For every $x\in\mathbb{R}^3$, $|P_Ux|^2+|P_Lx|^2=|x|^2$, and each of the projections has norm $1$  in the space  $L_c(\mathbb{R}^3,\mathbb{R}^3)$ of linear (continuous) maps $\mathbb{R}^3\to\mathbb{R}^3$.  

For a function $f:\mathbb{R}^n\supset dom\to\mathbb{R}^k$ on an open subset derivatives as linear maps $\mathbb{R}^n\to\mathbb{R}^k$ are denoted by $Df(x)$. For $n=k=1$,
$f'(x)=Df(x)1$. For partial derivatives in case $k=1$,  $\partial_jf(x)=Df(x)e_j$ for $j=1,\ldots,n$.

Let $M\subset{R}^n$ be a continuously differentiable submnifold.
For $x\in M$ the tangent space $T_xM$ is the set of tangent vectors $v=c'(0)$ of continuously differentiable curves $c:I\to\mathbb{R}^n$ with $I$ an interval, not a singleton, $c(I)\subset M$, $0\in I$, $c(0)=x$. A continuously differentiable map $f:M\supset dom\to N$, $dom$ open in $M$ and $N$ a continuously differentiable submnifold of $\mathbb{R}^k$, is locally given by restrictions of continuously differentiable maps $g:\mathbb{R}^n\supset U\to\mathbb{R}^k$. For such $U$ and $g$, and for $x\in dom\cap U$, the derivative of $f$ at $x$ is the linear map $T_xf:T_xM\to T_{f(x)}N$ 
given by $T_xf(v)=(g\circ c)'(0)=Dg(x)v$ for $v=c'(0)$ and $c$ and $g$ as above (with $c(I)\subset dom\cap U$).

The flow
${\mathcal F}$ generated by a vectorfield ${\mathcal V}:\mathbb{R}^n\supset{\mathcal U}\to\mathcal{R}^n$ which is locally Lipschitz continuous is the map $\mathbb{R}\times\mathbb{R}^n\supset dom_{{\mathcal F}}\to\mathbb{R}^n$ 
which is given by $(t,x)\in\,dom_{{\mathcal F}}$ if and only if $t$ belongs to the domain of the maximal solution $y:I_x\to\mathbb{R}^n$ of the differential equation $x'(t)={\mathcal V}(x(t))$ with initial value $y(0)=x$, and for such $y$, ${\mathcal F}(t,x)=y(t)$. ${\mathcal F}$ is of the same order of differentiability as ${\mathcal V}$. A subset $M\subset U$ is invariant under ${\mathcal F}$ if $x\in M$ implies ${\mathcal F}(t,x)\in M$ for all $t\in\mathbb{R}$ with $(t,x)\in \,dom_{{\mathcal F}}$. For a further subset $N\subset U$ the set $M$ is called invariant under ${{\mathcal F}}$ in $N$ if for every $x\in M\cap N$ and for every interval $I\ni 0$ with ${{\mathcal F}}(I\times\{x\})\subset N$ we have ${{\mathcal F}}(I\times\{x\})\subset M$.
A flowline $\xi:\mathbb{R}\to\mathbb{R}^n$ of a flow ${\mathcal F}$ on  $dom_{{\mathcal F}}=\mathbb{R}\times\mathbb{R}^n$ satisfies $\xi(t+s)={\mathcal F}(t,\xi(s))$ for all reals $s,t$.

\section{Transformation, scaling, projected flowlines}

In the appendix Section 9 we describe how the flow $F_V$ and the homoclinic flowline $h_V$ can be transformed to a continuously differentiable flow $F:\mathbb{R}\times\mathbb{R}^3\to\mathbb{R}^3$ and a flowline $h:\mathbb{R}\to\mathbb{R}^3$ of $F$ with the following properties. 
\begin{itemize}
\item[]\textbf{(F1)}
For all $t\in\mathbb{R}$, $F(t,0)=0$.
\item[]\textbf{(F2)}
$h(t)\neq0\neq h'(t)$ for all $t\in\mathbb{R}$, and $\lim_{|t|\to\infty}h(t)=0.$
\item[]\textbf{(F3)}
Every linear map $T(t):\mathbb{R}^3\ni x\mapsto D_2F(t,0)x\in\mathbb{R}^3$, $t\in\mathbb{R},$ satisfies  $T(t)L\subset L$ and $T(t)U\subset U$. For $t\in\mathbb{R},x\in\mathbb{R}^3,$ and $y=T(t)x$, 
\begin{equation*}
\left(
\begin{array}{c}
	y_1\\y_2
\end{array}
\right)
=e^{\sigma t}
\left(
\begin{array}{cc}
	\cos(\mu t) & \sin(\mu t)\\
	 -\sin(\mu t) & \cos(\mu t)
\end{array}
\right)
\left(
\begin{array}{c}
	x_1\\x_2
\end{array}
\right)\quad\mbox{and}\quad y_3=e^{ut}x_3.
\end{equation*}
\item[]\textbf{(F4)}
There exists $r_F>0$ such that $L$ and $U$ are invariant under $F$ in $\{x\in\mathbb{R}^3:|x|<r_F\}$.
\item[]\textbf{(F5)}
There exist reals $t_U<t_L$ with $h(t)\in U$ for all $t\le t_U$ and $h(t)\in L$ for all $t\ge t_L$.  Either $h(t)\in(0,\infty)e_3$ for all $t\le t_U$, or  $h(t)\in(-\infty,0)e_3$ for all $t\le t_U$.
\end{itemize}
In the sequel we focus on the case $h(t)\in(0,\infty)e_3$ for all $t\le t_U$, the other case being analogous.
We define
\begin{eqnarray}
B_1 & = & \{x\in\mathbb{R}^3:|P_Lx|\le1,|P_Ux|\le1\},\nonumber\\
 r_B & = & 2(\max_{0\le t\le1}|T(t)|+e^{u}+2),\quad\mbox{and}\nonumber\\
B & = & \{x\in\mathbb{R}^3:|x|\le r_B\}.\nonumber
\end{eqnarray}
Obviously, $B_1\subset B$. For $\epsilon>0$ we consider the scaled flows given by $F_{\epsilon}(t,x)=\frac{1}{\epsilon}F(t,\epsilon x)$. If
$x$ is a flowline of $F_{\epsilon}$ then $\epsilon x$ is a flowline of $F$, and conversely, if $y$ is a flowline of $F$, then $\frac{1}{\epsilon}y$ is a flowline of $F_{\epsilon}$. Observe that for all $\epsilon>0,x\in\mathbb{R}^3,y\in\mathbb{R}^3,t\in\mathbb{R}$,
$$
D_2F_{\epsilon}(t,x)y=\frac{1}{\epsilon}D_2F(t,\epsilon x)\epsilon y=D_2F(t,\epsilon x)y.
$$

\begin{proposition}
(i) For every $\eta>0$ there exists $\epsilon(\eta)>0$ such that for all $\epsilon\in(0,\epsilon(\eta))$,
$$
|D_2F_{\epsilon}(t,x)-T(t)|<\eta\,\,\mbox{for all}\,\,x\in B_1.
$$
(ii)
There exists $\epsilon_B>0$ such that $0<\epsilon<\epsilon_B$,
$$
F_{\epsilon}([0,1]\times B_1)\subset B,
$$
and $L$ and $U$ are invariant under $F_{\epsilon}$ in $B$.
\end{proposition}

{\bf Proof.}  1. On (i). By compactness and continuity there exists $\epsilon_(\eta)>0$ such that for $0<\epsilon<\epsilon(\eta)$ and $x\in B_1$ and for all $t\in[0,1]$,
$$
\eta\ge|D_2F(t,\epsilon x)-D_2F(t,0)|=|D_2F_{\epsilon}(t,x)-T(t)|.
$$

2. Let $0\le t\le1$. For $0<\epsilon<\epsilon(1)$ assertion (i) yields
$|D_2F_{\epsilon}(t,y)|\le\max_{0\le s\le1}|T(s)|+1\le r_B/2$ for all $y\in B_1$.  For $x\in B_1$ we have $|x|\le2$, and we infer
$$
|F_{\epsilon}(t,x)|=|F_{\epsilon}(t,x)-0|=|F_{\epsilon}(t,x)-F_{\epsilon}(t,0)|=\left|\int_0^1D_2F_{\epsilon}(t,sx)xds\right|\le r_B|x|/2\le r_B,
$$
hence $F_{\epsilon}(t,x)\in B$.

3. On invariance. Let $0<\epsilon<r_F/r_B$. Assume $x\in L$ and $F_{\epsilon}(s,x)\in B$ for $s$ between $0$ and $t\in\mathbb{R}$. Then we have $F(s,\epsilon x)\in \epsilon B$, or $|F(s,\epsilon x)|\le\epsilon r_B<r_F$. This yields $F(t,\epsilon x)\in L$. Consequently, $F_{\epsilon}(t,x)=\frac{1}{\epsilon}F(t,\epsilon x)\in\frac{1}{\epsilon}L=L$. Analogously for $x\in U$. $\Box$

\medskip

The next proposition expresses closeness of the flow to its linearization at the origin in terms of components in $U$ and $L$, for $0\le t\le 1$ and $x\in B_1$.

\begin{proposition}
For $\eta>0$, $0<\epsilon<\min\{\epsilon(\eta),\epsilon_B\}$, $x\in B_1$, and $0\le t\le1$, 
$$
|P_UF_{\epsilon}(t,x)-T(t)P_Ux)|\le\eta|P_Ux|\quad\mbox{and}\quad|P_LF_{\epsilon}(t,x)-T(t)P_Lx)|\le\eta|P_Lx|.
$$
Moreover,
$$
|P_UF_{\epsilon}(t,x)|\in (e^{ut}+[-\eta,\eta])|P_Ux|\quad\mbox{and}\quad|P_LF_{\epsilon}(t,x)|\in(e^{\sigma t}+[-\eta,\eta])|P_Lx|.
$$
\end{proposition}

{\bf Proof.} Let $x\in B_1$ and $0\le t\le1$. We apply Proposition 2.1(ii) to $P_Ux\in B_1$ and to $P_Lx\in B_1$ and obtain $F_{\epsilon}(t,P_Ux)\in U$ and $F_{\epsilon}(t,P_Lx)\in L$, hence $P_LF_{\epsilon}(t,P_Ux)=0$ and $P_UF_{\epsilon}(t,P_Lx)=0$. It follows that
\begin{eqnarray*}
|P_UF_{\epsilon}(t,x)-T(t)P_Ux| & = & |P_UF_{\epsilon}(t,x)-T(t)P_Ux-(P_UF_{\epsilon}(t,P_Lx)-P_UT(t)P_Lx)|\\
& = & \left|\int_0^1(P_UD_2F_{\epsilon}(t,P_Lx+s(x-P_Lx))[x-P_Lx]-P_UT(t)[x-P_Lx])ds\right|\\
& = & \left|\int_0^1P_U\{D_2F_{\epsilon}(t,P_Lx+s(x-P_Lx))[P_Ux]-P_UT(t)[P_Ux]\})ds\right|\\
& \le & |P_U|\max_{y\in B_1}|D_2F_{\epsilon}(t,y)-T(t)||P_Ux|\le\eta|P_Ux|\\
& & \mbox{(with Proposition 2.1(i) and}\quad|P_U|=1).
\end{eqnarray*}
We have $(P_Ux)_3=x_3$, and for $z=T(t)P_Ux\in U$, $|z|=|z_3|$ with $z_3=e^{ut}x_3$. Using the previous estimate we 
obtain
\begin{eqnarray*}
(e^{ut}-\eta)|x_3| & = & |z|-\eta|x_3|=|T(t)P_Ux|-\eta|P_Ux|\le|P_UF_{\epsilon}(t,x)|\\
& \le & |T(t)P_Ux|+\eta|P_Ux|=|z|+\eta|x_3|=(e^{ut}+\eta)|x_3|
\end{eqnarray*}
which yields
$$
|P_UF_{\epsilon}(t,x)|\in (e^{ut}+[-\eta,\eta])|P_Ux|.
$$
The remaining assertions are shown analogously. $\Box$

\begin{proposition}
Assume $0<\eta<e^{\sigma}$, $e^{\sigma}+\eta<1$, and $0<\epsilon<\min\{\epsilon(\eta),\epsilon_B\}$. Let $n\in\mathbb{N}$ and $x\in B_1$ be given with $|P_UF_{\epsilon}(j+1,x)|\le1$ for $j=0,\ldots,n-1$. Then we have $F_{\epsilon}(j,x)\in B_1$ for $j=0,\ldots,n$, and for $0\le t\le1$,
\begin{eqnarray*}
e^{\sigma(t+n)}(1-\eta e^{-\sigma})^{n+1}|P_Lx| & \le|P_LF_{\epsilon}(t+n,x)|\le & e^{\sigma(t+n)}(1+\eta e^{-\sigma})^{n+1}|P_Lx|,\\
e^{u(t+n)}(1-\eta)^{n+1}|P_Ux| & \le|P_UF_{\epsilon}(t+n,x)|\le & e^{u(t+n)}(1+\eta)^{n+1}|P_Ux|.
\end{eqnarray*}
\end{proposition}

{\bf Proof.} 1. From Proposition 2.2 for $t=1$, 
$$
 e^u(1-\eta)|P_Ux|\le(e^u-\eta)|P_Ux|\le|P_UF_{\epsilon}(1,x)|\le(e^u+\eta)|P_Ux|\le e^u(1+\eta)|P_Ux|
$$
and
$$
 e^{\sigma}(1-\eta e^{-\sigma})|P_Lx|=(e^{\sigma}-\eta)|P_Lx|\le|P_LF_{\epsilon}(1,x)|\le(e^{\sigma}+\eta)|P_Lx|= e^{\sigma}(1+\eta e^{-\sigma})|P_Lx|.
$$
Using $x\in B_1$ and $e^{\sigma}+\eta<1$ we infer $|P_LF_{\epsilon}(1,x)|\le 1$. With $|P_UF_{\epsilon}(1,x)|\le1$, we get $F_{\epsilon}(1,x)\in B_1$.

\medskip

By induction we obtain $F_{\epsilon}(j,x)\in B_1$ for $j=1,\ldots,n$, with
\begin{eqnarray*}
e^{\sigma j}(1-\eta e^{-\sigma})^ j|P_Lx| & \le|P_LF_{\epsilon}(j,x)|\le & e^{\sigma j}(1+\eta e^{-\sigma})^j|P_Lx|\quad\mbox{and}\\
e^{u j}(1-\eta)^j|P_Ux| & \le|P_UF_{\epsilon}(j,x)|\le & e^{u j}(1+\eta)^j|P_Ux|.
\end{eqnarray*}

2. On the lower estimates of the assertion.  For $j=n$ and $0\le t\le 1$ Proposition 2.2 yields
$$
(e^{u t}-\eta)|P_UF_{\epsilon}(n,x)|\le|P_UF_{\epsilon}(t+n,x)|.
$$
Use $ e^{u t}(1-\eta)\le e^{u t}-\eta$ and the lower estimate for $|P_UF_{\epsilon}(n,x)|$ in order to get 
$$
e^{u(t+n)}(1-\eta)^{n+1}|P_Ux|\le|P_UF_{\epsilon}(t+n,x)|.
$$
In the same way, now using $e^{\sigma t}(1-\eta e^{-\sigma})\le e^{\sigma t}-\eta$, one finds 
$$
e^{\sigma(t+n)}(1-\eta e^{-\sigma})^{n+1}|P_Lx|\le|P_LF_{\epsilon}(t+n,x)|.
$$
3. The upper estimates of the assertion are shown analogously. $\Box$

\medskip

For later use we turn to exponential estimates for flowlines which stay sufficiently long in $B_1$.  

\begin{proposition}
Let $\tilde{\eta}>0$ be given. Assume  $0<\eta<e^{\sigma}$ and $e^{\sigma}+\eta<1$ as in Proposition 2.3, and in addition
$$
\log(1+\eta e^{-\sigma})<\tilde{\eta}\quad\mbox{and}\quad\log\left(\frac{1}{1-\eta e^{-\sigma}}\right)<\tilde{\eta}.
$$
Let $n\in\mathbb{N}$ be given with
$$
\frac{n+1}{n}\log(1+\eta e^{-\sigma})<\tilde{\eta}\quad\mbox{and}\quad\frac{n+1}{n}\log\left(\frac{1}{1-\eta e^{-\sigma}}\right)<\tilde{\eta}.
$$
Let $0<\epsilon<\min\{\epsilon(\eta),\epsilon_B\}$ and consider
$x\in B_1$ with $|P_UF_{\epsilon}(j+1,x)|\le1$ for $j=0,\ldots,n-1$ as in Proposition 2.3. Then we have, for $0\le t\le1$,
\begin{eqnarray*}
e^{(\sigma-\tilde{\eta})(t+n)}|P_Lx| & \le|P_LF_{\epsilon}(t+n,x)|\le & e^{(\sigma+\tilde{\eta})(t+n)}|P_Lx|,\\
e^{(u-\tilde{\eta})(t+n)}|P_Ux| & \le|P_UF_{\epsilon}(t+n,x)|\le & e^{(u+\tilde{\eta})(t+n)}|P_Ux|.
\end{eqnarray*}
\end{proposition}

{\bf Proof.} In view of Proposition 2.3 the estimates of $|P_UF_{\epsilon}(t+n,x)|$ follow from the estimates
$$
 e^{(u-\tilde{\eta})(t+n)}\le e^{u(t+n)}(1-\eta)^{n+1}\quad\mbox{and}\quad e^{u(t+n)}(1+\eta)^{n+1}\le  e^{(u+\tilde{\eta})(t+n)},
$$
or equivalently, 
$$
-(t+n)\tilde{\eta}\le(n+1)\log(1-\eta)\quad\mbox{and}\quad(n+1)\log(1+\eta)\le(t+ n)\tilde{\eta}.
$$
Sufficient for the latter are
$$
-n\tilde{\eta}\le(n+1)\log(1-\eta)\quad\mbox{and}\quad(n+1)\log(1+\eta)\le n\tilde{\eta},
$$
which are a consequence of the hypotheses on $\eta$ and $n$  in combination with
$$
\log\left(\frac{1}{1-\eta}\right)<\log\left(\frac{1}{1-\eta e^{-\sigma}}\right)\quad\mbox{and}\quad
\log(1+\eta)<\log(1+\eta e^{-\sigma}).
$$

Similarly the estimates of $|P_LF_{\epsilon}(t+n,x)|$ follow from the estimates
$$
e^{(\sigma-\tilde{\eta})(t+n)}\le e^{\sigma(t+n)}(1-\eta  e^{-\sigma})^{n+1}\quad\mbox{and}\quad e^{\sigma(t+n)}(1+\eta e^{-\sigma})^{n+1}\le e^{(\sigma+\tilde{\eta})(t+n)}
$$
which are equivalent to
$$
-\tilde{\eta}(t+n)\le(n+1)\log(1-\eta e^{-\sigma})\quad\mbox{and}\quad (n+1)\log(1+\eta e^{-\sigma})\le\tilde{\eta}(t+n).
$$
Sufficient for the latter are 
$$
-\tilde{\eta}n\le(n+1)\log(1-\eta e^{-\sigma})\quad\mbox{and}\quad (n+1)\log(1+\eta e^{-\sigma})\le\tilde{\eta}n.
$$
which are obvious from the hypotheses on $\eta$ and $n$. $\Box$

\section{Angles}

This section deals with angles in the plane $L$, along projected flowlines $t\mapsto P_LF_{\epsilon}(t,x)$, for $x\in B_1\setminus U$. The first step computes such angles for $0\le t\le1$.

\begin{proposition}
Assume $0<\eta<\frac{e^{\sigma}}{2}$, $0<\epsilon<\min\{\epsilon(\eta),\epsilon_B\}$, $x\in B_1\setminus U$, and
$$
\frac{1}{|P_Lx|}P_Lx=\left(\begin{array}{c}\cos(\psi)\\ \sin(\psi)\\0\end{array}\right)\quad\mbox{for some}\quad\psi\in\mathbb{R}.
$$
Let $0\le t\le 1$. For the unique $\Delta=\Delta(t,x,\psi)\in[-\pi,\pi)$ with
$$
\frac{1}{|P_LF_{\epsilon}(t,x)|}P_LF_{\epsilon}(t,x)=\left(\begin{array}{c}\cos(\psi-\mu t+\Delta)\\ \sin(\psi-\mu t+\Delta)\\0\end{array}\right)
$$
we have
$$
\Delta(t,x,\psi)=\arcsin(<v,w^{\perp}>)\in\left(-\frac{\pi}{2},\frac{\pi}{2}\right)
$$
where
$$
v=\frac{1}{|P_LF_{\epsilon}(t,x)|}P_LF_{\epsilon}(t,x),\quad w=\frac{1}{|T(t)P_Lx|}T(t)P_Lx,\quad\mbox{and}\quad w^{\perp}=\left(\begin{array}{c}-w_2\\w_1\\0\end{array}\right).
$$
In particular, $\Delta(0,x,\psi)=0$. Moreover,
$$
|<v,w^{\perp}>|\le\frac{\eta}{e^{\sigma}-\eta}\quad\mbox{and}\quad|\Delta|\le\arcsin\left(\frac{\eta}{e^{\sigma}-\eta}\right).
$$
\end{proposition}

Compare Figure 1, bottom.

\medskip

{\bf Proof.} 1. Using $F_{\epsilon}(B_1\times\{x\})\subset B$ and invariance of $U$ under $F_{\epsilon}$ in $B$ from  Proposition 2.1 (ii) we infer $F_{\epsilon}(t,x)\in B\setminus U$, hence $P_LF_{\epsilon}(t,x)\neq0$. Also, $T(t)P_Lx\neq0$.

2. We have
\begin{eqnarray*}
\left(\begin{array}{c}w_1\\w_2\end{array}\right) & = & \left(\begin{array}{cc}\cos(\mu t) & \sin(\mu t)\\-\sin(\mu t) & \cos(\mu t)\end{array}\right)\cdot\left(\begin{array}{c}\cos(\psi)\\ \sin(\psi)\end{array}\right)=\left(\begin{array}{c}\cos(\psi-\mu t)\\ \sin(\psi-\mu t)\end{array}\right)\\
& = & \left(\begin{array}{cc}\cos(\psi-\mu t) & -\sin(\psi-\mu t)\\ \sin(\psi-\mu t) & \cos(\psi-\mu t)\end{array}\right)\cdot\left(\begin{array}{c}1\\0\end{array}\right),
\end{eqnarray*}
hence
$$
\left(\begin{array}{c}-w_2\\w_1\end{array}\right)=\left(\begin{array}{cc}\cos(\psi-\mu t) & -\sin(\psi-\mu t)\\ \sin(\psi-\mu t) & \cos(\psi-\mu t)\end{array}\right)\cdot\left(\begin{array}{c}0\\1\end{array}\right).
$$
It follows that multiplication with the matrix
$$
\left(\begin{array}{ccc}\cos(\psi-\mu t) & -\sin(\psi-\mu t) & 0\\ \sin(\psi-\mu t) & \cos(\psi-\mu t) & 0\\0 & 0 & 1\end{array}\right)^{-1}=\left(\begin{array}{ccc}\cos(\psi-\mu t) & \sin(\psi-\mu t) & 0\\-\sin(\psi-\mu t) & \cos(\psi-\mu t) & 0\\0 & 0 & 1\end{array}\right)
$$
defines a linear map $\rho:L\to L$ which satisfies $\rho w=e_1$ and $\rho w^{\perp}=e_2$. The map $\rho$ preserves the inner product and the norm.
We obtain
\begin{eqnarray*}
\left(\begin{array}{c}\cos\\ \sin\\0\end{array}\right)(\Delta) & = & \left(\begin{array}{c}\cos(\psi-\mu t+\Delta-(\psi-\mu t))\\ \sin(\psi-\mu t+\Delta-(\psi-\mu t))\\0\end{array}\right)\\
 & = & \left(\begin{array}{ccc}\cos(\psi-\mu t) & \sin(\psi-\mu t) & 0\\ -\sin(\psi-\mu t) & \cos(\psi-\mu t) & 0\\ 0 & 0 & 1\end{array}\right)\cdot\left(\begin{array}{c}\cos(\psi-\mu t+\Delta)\\ \sin(\psi-\mu t+\Delta)\\0\end{array}\right)\\
& = & \rho v=<\rho v,e_1>e_1+<\rho v,e_2>e_2=<\rho v,\rho w>e_1+<\rho v,\rho w^{\perp}>e_2\nonumber\\
& = & <v,w>e_1+<v,w^{\perp}>e_2,
\end{eqnarray*}
hence $\cos(\Delta)=<v,w>$ and $\sin(\Delta)=<v,w^{\perp}>$. Let $\tilde{v}=P_LF_{\epsilon}(t,x)\neq0$ and $\tilde{w}=T(t)P_Lx\neq0$. Then
$$
<v,w>=\frac{1}{|\tilde{v}||\tilde{w}|}<\tilde{v},\tilde{w}>.
$$
Proposition 2.2 shows that $\tilde{z}=\tilde{v}-\tilde{w}$ satisfies $|\tilde{z}|\le\eta|P_Lx|$. Recall $|\tilde{w}|=e^{\sigma t}|P_Lx|$. We infer
$$
<\tilde{v},\tilde{w}>=<\tilde{w}+\tilde{z},\tilde{w}>\,\,\ge|\tilde{w}|^2-|\tilde{z}||\tilde{w}|\ge|P_Lx|^2e^{\sigma t}(e^{\sigma t}-\eta)\ge|P_Lx|^2e^{\sigma}(e^{\sigma}-\eta)>0,
$$
which yields $\cos(\Delta)=<v,w>>0$. Recall $-\pi\le\Delta<\pi$. It follows that $|\Delta|<\frac{\pi}{2}$. Consequently,
$$
\Delta=\arcsin(\sin(\Delta))=\arcsin(<v,w^{\perp}>).
$$
The equation $\Delta(0,x,\psi)=0$ follows from $P_LF_{\epsilon}(0,x)=P_Lx=T(0)P_Lx$, which yields $v=w$, hence $<v,w^{\perp}>=0$.

3. Proof of the estimates of $<v,w^{\perp}>$ and $\Delta$. We have 
$$
|<v,w^{\perp}>|=\frac{1}{|\tilde{v}||\tilde{w}|}|<\tilde{v},\tilde{w}^{\perp}>|=\frac{1}{|\tilde{v}||\tilde{w}|}|<\tilde{z},\tilde{w}^{\perp}>|\le\frac{|\tilde{z}|}{|\tilde{v}|}.
$$
Using $|\tilde{w}|=e^{\sigma t}|P_Lx|$ and $|\tilde{z}|\le\eta|P_Lx|$ we get
$$
|\tilde{v}|\ge|\tilde{w}|-|\tilde{z}|=e^{\sigma t}|P_Lx|-|\tilde {z}|\ge(e^{\sigma t}-\eta)|P_Lx|\ge(e^{\sigma}-\eta)|P_Lx|
$$
and obtain
$$
\frac{|\tilde{z}|}{|\tilde{v}|}\le\frac{\eta}{e^{\sigma}-\eta}<1
$$
where the last inequality holds by the hypothesis on $\eta$. Finally,
$$
|\Delta|=|\arcsin(<v,w^{\perp}>)|=\arcsin(|<v,w^{\perp}>|)\le\arcsin\left(\frac{\eta}{e^{\sigma}-\eta}\right).\quad\Box
$$

The next result describes the angles along projected flowlines $P_LF_{\epsilon}(\cdot,x)$, $x\in B_1\setminus U$, by continuous functions, a little longer than the projected values $P_UF_{\epsilon}(\nu,x)$, $\nu\in\mathbb{N}$, remain in $B_1\setminus U$. For later use we restrict attention to flowlines which start from the strip on the cylinder
$$
M_I=\{x\in\mathbb{R}^3:|P_Lx|=1\}
$$
which is given by $0<x_3\le1$, 
\begin{equation}
x=\left(\begin{array}{c}\cos(\psi)\\ \sin(\psi)\\ \delta\end{array}\right)\quad\mbox{with}\quad\psi\in\mathbb{R}\quad\mbox{and}\quad0<\delta\le1.
\end{equation}
See Figure 2, top.

\medskip

\begin{figure}
	\includegraphics[page=1,scale=0.7]{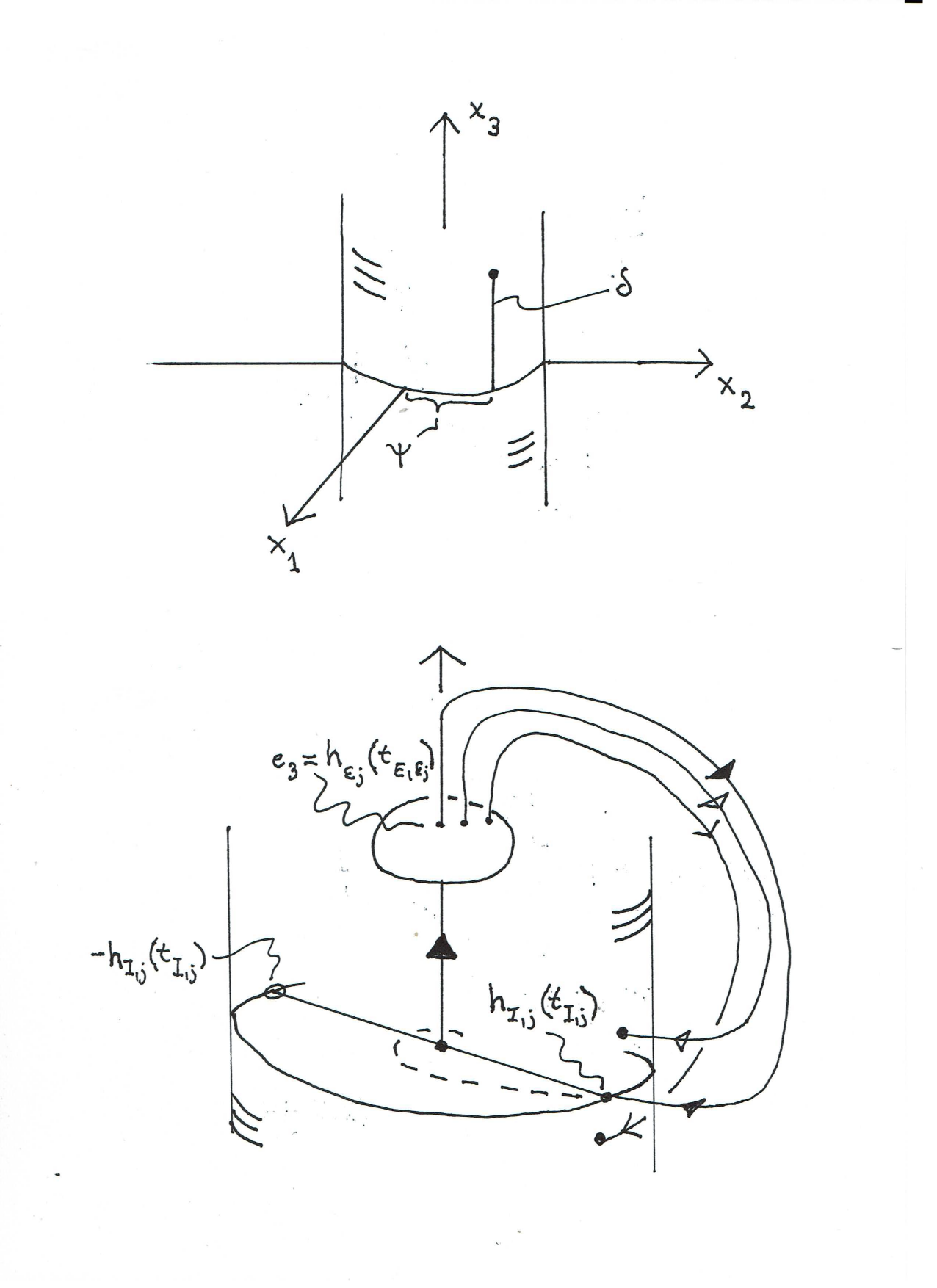}
	\caption{Top: The relations (2). Bottom: The exterior map.}
\end{figure}

For integers $n\ge0$ we define
$$
dom_n=\{(t,\psi,\delta)\in[0,n+1]\times\mathbb{R}\times(0,1]: |P_UF_{\epsilon}(\nu,x)|\le1\,\,\mbox{for}\,\,
x\,\,\mbox{given by }\,\,(2)\,\,\mbox{and}\,\,\nu=0,\ldots,n\}.
$$

\begin{proposition}
Assume $0<\eta<\frac{e^{\sigma}}{2}$, $0<\epsilon<\min\{\epsilon(\eta),\epsilon_B\}$.

\medskip

(i) For every integer $n\ge0$ there exists a continuous function $\phi^{(n)}:dom_n\to\mathbb{R}$ so that for each $(t,\psi,\delta)\in dom_n$,
with $x$ given by (2), we have
\begin{equation}
\frac{1}{|P_LF_{\epsilon}(t,x)|}P_LF_{\epsilon}(t,x)=\left(\begin{array}{c}\cos(\phi^{(n)}(t,\psi,\delta))\\ \sin(\phi^{(n)}(t,\psi,\delta))
\\0\end{array}\right),
\end{equation}
and in case $n\le t\le n+1$,
\begin{equation}
\psi-t\mu-(n+1)\arcsin\left(\frac{\eta}{e^{\sigma}-\eta}\right)\le\phi^{(n)}(t,\psi,\delta)\le\psi-t\mu+(n+1)\arcsin\left(\frac{\eta}{e^{\sigma}-\eta}\right).
\end{equation}

(ii) For every $n\in\mathbb{N}$, $([0,n]\times\mathbb{R}\times(0,1])\cap dom_n\subset dom_{n-1}$, and on this set,
$\phi^{(n)}(t,\psi,\delta)=\phi^{(n-1)}(t,\psi,\delta)$.
\end{proposition}

{\bf Proof.} 1. Obviously, $([0,n]\times\mathbb{R}\times(0,1])\cap dom_n\subset dom_{n-1}$ for all integers $n\ge1$.

\medskip

2. Proof of assertion (i). We construct the functions $\phi^{(n)}$ recursively. 

\medskip

2.1. For $n=0$ Proposition 3.1 shows that the continuous function $\phi^{(0)}:dom_0\to\mathbb{R}$ defined by
$$
\phi^{(0)}(t,\psi,\delta)=\psi-t\mu+\arcsin(\Delta)\quad\mbox{for}\quad0\le t\le1,
$$
with $\Delta=\Delta(t,x,\psi)$ from Proposition 3.1, for $x$ given by (2), satisfies (3) and (4). 

\medskip

2.2. Suppose now that for some integer $n\ge0$ the continuous function $\phi^{(n)}:dom_n\to\mathbb{R}$ satisfies (3) and (4). We define $\phi^{(n+1)}:dom_{n+1}\to\mathbb{R}$ as follows. 

\medskip

For $(t,\psi,\delta)\in dom_{n+1}$ with $0\le t\le n+1$ we have $(t,\psi,\delta)\in dom_n$ due to Part 1, and we  set $\phi^{(n+1)}(t,\psi,\delta)=\phi^{(n)}(t,\psi,\delta)$.  

\medskip

For $(t,\psi,\delta)\in dom_{n+1}$ with $n+1\le t\le n+2$, observe first that for $x$ given by (2) the property $|P_UF_{\epsilon}(\nu,x)|\le1$ for $\nu=0,\ldots,n+1$ yields
$F_{\epsilon}(\nu,x)\in B_1$ for $\nu=1,\ldots,n+1$, by means of Proposition 2.3. By Proposition 2.1 (ii), $F_{\epsilon}(s,x)\in B$ for $0\le s\le n+2$. Using this in combination with $x\in B_1\setminus U\subset B\setminus U$  we get $F_{\epsilon}(n+1,x)\in B\setminus U$, by invariance of $U$ from Proposition 2.1 (ii). Altogether, 
$y=F_{\epsilon}(n+1,x)$ is contained in $B_1\setminus U$, and
$$
\frac{1}{|P_Ly|}P_Ly=\left(\begin{array}{c}\cos(\phi^{(n)}(n+1,\psi,\delta))\\ \sin(\phi^{(n)}(n+1,\psi,\delta))\\0\end{array}\right).
$$
An application of Proposition 3.1 to $y\in B_1\setminus U$ shows that the continuous function 
$$
\phi^{\ast}:dom_0\to\mathbb{R},\quad\phi^{\ast}(t,\psi,\delta)=\phi^{(n)}(n+1,\psi,\delta)-t\mu+\arcsin(\Delta)\quad\mbox{for}\quad0\le t\le1,
$$
with $\Delta=\Delta(t,y,\phi^{(n)}(n+1,\psi,\delta))$ according to Proposition 3.1, satisfies (3) and (4) with $y$ and $\phi^{(n)}(n+1,\psi,\delta)$  in place of $x$ and $\psi$, respectively. Also,
$$
\phi^{(n)}(n+1,\psi,\delta)-t\mu-\arcsin\left(\frac{\eta}{e^{\sigma}-\eta}\right)\le\phi^{\ast}(t,\psi,\delta)\le\phi^{(n)}(n+1,\psi,\delta)-t\mu+\arcsin\left(\frac{\eta}{e^{\sigma}-\eta}\right).
$$

For $(t,\psi,\delta)\in dom_{n+1}$ with $n+1\le t\le n+2$ we have $(t-(n+1),\psi,\delta)\in dom_0$. We complete the definition of $\phi^{(n+1)}$ by
$$
\phi^{(n+1)}(t,\psi,\delta)=\phi^{\ast}(t-(n+1),\psi,\delta)\quad\mbox{for}\quad(t,\psi,\delta)\in dom_{n+1}\quad\mbox{with}\quad n+1<t\le n+2.
$$
Proof that $\phi^{(n+1)}$ is continuous. Let $(t,\psi,\delta)\in dom_{n+1}$ and a sequence $(t_m,\psi_m,\delta_m)_{m\in\mathbb{N}}$ in $dom_{n+1}$ with $\lim_{m\to\infty}(t_m,\psi_m,\delta_m)=(t,\psi,\delta)$  be given. It is enough to show that for a subsequence we have $\phi^{(n+1)}(t_{m_{\mu}},\psi_{m_{\mu}},\delta_{m_{\mu}})\to\phi^{(n+1)}(t,\psi,\delta)$ as $\mu\to\infty$. In case $t<n+1$, $(t,\psi,\delta)\in dom_n$ and $\phi^{(n+1)}(t,\psi,\delta)=\phi^{(n)}(t,\psi,\delta)$,
and for $m$ sufficiently large, $t_m<n+1$, hence $(t_m,\psi_m,\delta_m)\in dom_n$ and  $\phi^{(n+1)}(t_m,\psi_m,\delta_m)=\phi^{(n)}(t_m,\psi_m,\delta_m)$ for such $m$, and continuity of $\phi^{(n)}$ yields\\
$\phi^{(n+1)}(t_m,\psi_m,\delta_m)\to\phi^{(n+1)}(t,\psi,\delta)$
for $m\to\infty$.
In case $n+1<t\le n+2$ a similar argument with $\phi^{\ast}(\cdot-(n+1),\cdot,\cdot)$ in place of $\phi^{(n)}$ yields $\phi^{(n+1)}(t_m,\psi_m,\delta_m)\to\phi^{(n+1)}(t,\psi,\delta)$ for $m\to\infty$ as well.
In case $t=n+1$ we distinguish the subcases that the set of indices $m$ with $t_m\le n+1$ is bounded or unbounded. If it is unbounded then we can argue for a subsequence as in case $t<n+1$ above. If it is bounded we can argue  as in case $n+1<t\le n+2$ above and find
\begin{eqnarray}
lim_{m\to\infty}\phi^{(n+1)}(t_m,\psi_m,\delta_m)&=&lim_{m\to\infty}\phi^{\ast}(t_m-(n+1),\psi_m,\delta_m)=\phi^{\ast}(0,\psi,\delta)\nonumber\\
& = &\phi^{(n)}(n+1,\psi,\delta)=\phi^{(n+1)}(n+1,\psi,\delta).\nonumber
\end{eqnarray}

2.3.  Proof of (3) for $\phi^{(n+1)}$. By the properties of $\phi^{(n)}$,  for $(t,\psi,\delta)\in dom_{n+1}$ with $0\le t\le n+1$, obviously 
$$
\frac{1}{|P_LF_{\epsilon}(t,x)|}P_LF_{\epsilon}(t,x)=\left(\begin{array}{c}\cos(\phi^{(n)}(t,\psi,\delta))\\ \sin(\phi^{(n)}(t,\psi,\delta))\\0\end{array}\right)=\left(\begin{array}{c}\cos(\phi^{(n+1)}(t,\psi,\delta))\\ \sin(\phi^{(n+1)}(t,\psi,\delta))\\0\end{array}\right).
$$
For $(t,\psi,\delta)\in dom_{n+1}$ with $n+1<t\le n+2$ let $s=t-(n+1)$ and $y=F_{\epsilon}(n+1,x)$ with $x$ given by (2). Using the version of (3) for $\phi^{\ast}$ we have
\begin{eqnarray*}
\frac{1}{|P_LF_{\epsilon}(t,x)|}P_LF_{\epsilon}(t,x) &  = & \frac{1}{|P_LF_{\epsilon}(s,y)|}P_LF_{\epsilon}(s,y)
=\left(\begin{array}{c}\cos(\phi^{\ast}(s,\psi,\delta))\\ \sin(\phi^{\ast}(s,\psi,\delta))\\0\end{array}\right)\\
& = & \left(\begin{array}{c}\cos(\phi^{\ast}(t-(n+1),\psi,\delta))\\ \sin(\phi^{\ast}(t-(n+1),\psi,\delta))\\0\end{array}\right)
=\left(\begin{array}{c}\cos(\phi^{(n+1)}(t,\psi,\delta))\\ \sin(\phi^{(n+1)}(t,\psi,\delta))\\0\end{array}\right).
\end{eqnarray*}

2.4.  Proof of (4) for $\phi^{(n+1)}$. For $(n+1+t,\psi,\delta)\in dom_{n+1}$ with $0<t\le1$ we have
\begin{eqnarray*}
\phi^{(n+1)}(n+1+t,\psi,\delta) & = & \phi^{\ast}(t,\psi,\delta)\le\phi^{(n)}(n+1,\psi,\delta)-t\mu+\arcsin\left(\frac{\eta}{e^{\sigma}-\eta}\right)\\
& \le & \left(\psi-(n+1)\mu+(n+1)\arcsin\left(\frac{\eta}{e^{\sigma}-\eta}\right)\right)-t\mu +\arcsin\left(\frac{\eta}{e^{\sigma}-\eta}\right)\\
& & \mbox{(by (4) for $\phi^{(n)}$)}\\
& = & \psi-(n+1+t)\mu+(n+2)\arcsin\left(\frac{\eta}{e^{\sigma}-\eta}\right)\nonumber.
\end{eqnarray*}
Analogously we get the lower estimate
$$
\psi-(n+1+t)\mu-(n+2)\arcsin\left(\frac{\eta}{e^{\sigma}-\eta}\right)\le\phi^{(n+1)}(n+1+t,\psi,\delta)
$$. 
for $0<t\le1$. By continuity both estimates hold also at $t=n+1$.

\medskip

2.5. The proof of assertion (i) is complete. Assertion (ii) is obvious from the previous construction which includes the relation
$$
\phi^{(n)}(t,\psi,\delta)=\phi^{(n-1)}(t,\psi,\delta)
$$  
for $(t,\psi,\delta)\in([0,n]\times\mathbb{R}\times(0,1])\cap dom_n\subset dom_{n-1}$. $\Box$

\medskip

We proceed to estimates of $\phi^{(n)}$ from Proposition 3.2 on the interval $[n,n+1]$ in terms of affine linear maps with slopes close to $\pm\mu$.

\begin{corollary}
Let $\tilde{\eta}>0$ and  $\eta>0$ be given with $\eta<\frac{e^{\sigma}}{2}$ and
$$
\arcsin\left(\frac{\eta}{e^{\sigma}-\eta}\right)<\tilde{\eta}.
$$
Let $0<\epsilon<\min\{\epsilon(\eta),\epsilon_B\}$ and let $n\in\mathbb{N}$ be given with
$$
\frac{n+1}{n}\arcsin\left(\frac{\eta}{e^{\sigma}-\eta}\right)<\tilde{\eta}.
$$
Then the function $\phi^{(n)}$ obtained in Proposition 3.2 satisfies
$$
\psi-(n+t)(\mu+\tilde{\eta})\le\phi^{(n)}(n+t,\psi,\delta)\le\psi-(n+t)(\mu-\tilde{\eta})
$$
for every $(n+t,\psi,\delta)\in dom_n$ with $0\le t\le1$.
\end{corollary}

{\bf Proof.} 1. The estimate (4) in Proposition 3.2 shows that the upper estimate of the assertion follows from 
$$
-(n+t)\mu+(n+1)\arcsin\left(\frac{\eta}{e^{\sigma}-\eta}\right)\le-(n+t)(\mu-\tilde{\eta})
$$
which is equivalent to
$$
(n+1)\arcsin\left(\frac{\eta}{e^{\sigma}-\eta}\right)\le(n+t)\tilde{\eta}.
$$
The preceding estimate follows from
$$
(n+1)\arcsin\left(\frac{\eta}{e^{\sigma}-\eta}\right)\le n\tilde{\eta}
$$
which is obvious from the hypotheses on $\eta$ and $n$.

2. Analogously the estimate (4) shows that the lower estimate of the assertion follows from
$$
-(n+t)\mu-(n+1)\arcsin\left(\frac{\eta}{e^{\sigma}-\eta}\right)\ge-(n+t)(\mu+\tilde{\eta})
$$
which is equivalent to
$$
-(n+1)\arcsin\left(\frac{\eta}{e^{\sigma}-\eta}\right)\ge-(n+t)\tilde{\eta}.
$$
The preceding estimate follows from
$$
-(n+1)\arcsin\left(\frac{\eta}{e^{\sigma}-\eta}\right)\ge- n\tilde{\eta}
$$
which is obvious from the hypotheses on $\eta$ and $n$. $\Box$

\section{Transversality, and exterior maps}

For every $\epsilon>0$ the flowline $h_{\epsilon}=\frac{1}{\epsilon}h$ of $F_{\epsilon}$ is homoclinic with
$h_{\epsilon}(t)\neq0\neq h_{\epsilon}'(t)$ everywhere and $h_{\epsilon}(t)\to0$ as $|t|\to\infty$.
 
\begin{proposition}
(i) For $0<\epsilon<|h(t_U)|$ there are reals $t_{E,\epsilon}\le t_U$ with 
$$
h_{\epsilon}(t_{E,\epsilon})=e_3\quad\mbox{and}\quad h_{\epsilon}(t)\in(0,\infty)e_3\quad\mbox{for all}\quad t\le t_{E,\epsilon}.
$$
(ii) There are strictly montonic sequences $(t_{I,j})_{j\in\mathbb{N}}$ in $[t_L,\infty)$ and $(\epsilon_j)_{j\in\mathbb{N}}$ in $(0,\infty)$ with $t_{I,j}\to\infty$ and $\epsilon_j\to0$ as $j\to\infty$ such that for every $j\in\mathbb{N}$,
$$
\epsilon_j<|h(t_U)|,\quad |h_{\epsilon_j}(t_{I,j})|=1,\quad (|h_{\epsilon_j}|^2)'(t_{I,j})<0,\quad h_{\epsilon_j}(t)\in L\quad\mbox{for}\quad t\ge t_{I,j},\quad h_{\epsilon_j}'(t_{I,j})\in L.
$$
\end{proposition}

{\bf Proof.} 1. On (i). Recall that for all $t\le t_U$, $h(t)\in(0,\infty)e_3$. Let $0<\epsilon<|h(t_U)|$.
The relation  $\lim_{t\to-\infty}h(t)=0$ shows that for each $\epsilon\in(0,|h(t_U)|)$ there exists $t_{E,\epsilon}\le t_U$ with $|h(t_{E,\epsilon})|=\epsilon$. It follows that $h(t_{E,\epsilon})=\epsilon e_3$, hence $h_{\epsilon}(t_{E,\epsilon})=e_3$. For $0<\epsilon<|h(t_U)|$ and for all $t\le t_{E,\epsilon}\le t_U$
we get $h_{\epsilon}(t)=\frac{1}{\epsilon}h(t)\in\frac{1}{\epsilon}(0,\infty)e_3=(0,\infty)\epsilon$.

\medskip

2. On (ii). From $h(t)\neq0$ everywhere and $h(t)\to0$ as $t\to\infty$ we get a strictly increasing sequence   $(t_{I,j})_{j\in\mathbb{N}}$ in $[t_L,\infty)$ with  $t_{I,j}\to\infty$ for $j\to\infty$ so that $|h(t_{I,j})|<|h(t_U)|$ and $(|h|^2)'(t_{I,j})<0$ for all $j\in\mathbb{N}$, and the sequence given by $\epsilon_j=|h(t_{I,j})|$ is strictly decreasing. It follows that 
$$
|h_{\epsilon_j}(t_{I,j})|=\frac{1}{\epsilon_j}|h(t_{I,j})|=1\quad\mbox{and}\quad
(|h_{\epsilon_j}|^2)'(t_{I,j})=\frac{1}{\epsilon_j^2}|(|h|^2)'(t_{I,j})<0\quad\mbox{for all}\quad j\in\mathbb{N}.
$$
Also, for $t\ge t_{I,j}\ge t_L$, $h_{\epsilon_j}(t)=\frac{1}{\epsilon_j}h(t)\in\frac{1}{\epsilon_j}L=L$, which yields
$h_{\epsilon_j}'(t_{I,j})\in L$. $\Box$

\medskip

 We want to describe the behaviour of flowlines close to the homoclinic loop $h_{\epsilon}(\mathbb{R})\cup\{0\}$  for small $\epsilon>0$. This will be done in terms of a return map which is given by the return of flowlines $F_{\epsilon}(\cdot,x)$ from points $x$ in the
cylinder
$$
M_I=\{z\in\mathbb{R}^3:|P_Lz|=1\}
$$ 
with $0<x_3$ slightly above the plane $L$, to targets in  $M_I$. The return map will be obtained as a composition of an {\it inner map}, which follows the flow until it reaches the plane
$$
M_E=e_3+L,
$$
parallel to $L$
% (compare Figure 1, top, on page ??)
, with an {\it exterior map} following the flow from a neighbourhood of $e_3$ in $M_E$ until it reaches $M_I$. 

\medskip

The constructions and continuous differentiability of the inner and exterior maps requires that the homoclinic flowline intersects the smooth 2-dimensional submanifolds $M_E$ and $M_I$ of $\mathbb{R}^3$ transversally. Obviously, 
$T_xM_E=L$ for all $x\in M_E$, and at $x\in M_I\cap L$, 
$$
T_xM_I=\mathbb{R}x^{\perp}\oplus\mathbb{R}e_3,\quad\mbox{with}\quad x^{\perp}=\left(\begin{array}{c}-x_2\\x_1\\0\end{array}\right).
$$. 

\begin{proposition}
(i) For $0<\epsilon<|h(t_U)|$,  $\partial_1F_{\epsilon}(0,e_3)=D_1F_{\epsilon}(0,e_3)1=h_{\epsilon}'(t_{E,\epsilon})\notin L=T_{e_3}M_E$.

\medskip
(ii) Let $j\in\mathbb{N}$ be given and  $x=h_{\epsilon_j}(t_{I,j})$. Then
$$
\partial_1F_{\epsilon_j}(0,x)=h_{\epsilon_j}'(t_{I,j})=\partial_1F_{\epsilon_j}(t_{I,j}-t_{E,\epsilon_j},e_3),\quad
<\partial_1F_{\epsilon_j}(0,x),x>\,\,<0,\quad\mbox{and}\quad \partial_1F_{\epsilon_j}(0,x)\notin T_xM_I.
$$
\end{proposition}

{\bf Proof.} 1. On (i).  From 
$$
F_{\epsilon}(t,e_3)-F_{\epsilon}(0,e_3)=F_{\epsilon}(t,h_{\epsilon}(t_{E,\epsilon}))-e_3=h_{\epsilon}(t+t_{E,\epsilon})-h_{\epsilon}(E,t_{\epsilon})
$$ 
we have $\partial_1F_{\epsilon}(0,e_3)=h_{\epsilon}'(t_{E,\epsilon})\neq0$. Using $h_{\epsilon}(t)\in U$ for $t\le t_ {E,\epsilon}$ we get $h_{\epsilon}'(t_{E,\epsilon})\in U$. It follows that $\partial_1F_{\epsilon}(0,e_3)\in U\setminus\{0\}\subset\mathbb{R}^3\setminus L$. 

\medskip

2. On (ii). Let $j\in\mathbb{N}$ be given. Let   $x=h_{\epsilon_j}(t_{I,j})\in M_I\cap L$. Arguing as in Part 1 we get $\partial_1F_{\epsilon_j}(0,x)=h_{\epsilon_j}'(t_{I,j})=\partial_1F_{\epsilon_j}(t_{I,j}-t_{E,\epsilon_j},e_3)$. From Proposition 4.1 (ii),
$$
0> (|h_{\epsilon_j}|^2)'(t_{I,j})=2<h_{\epsilon_j}'(t_{I,j}),h_{\epsilon_j}(t_{I,j})>=2<\partial_1F_{\epsilon_j}(0,x),x>
$$
and $h_{\epsilon_j}'(t_{I,j})\in L$.  From the preceding relations, $\partial_1F_{\epsilon_j}(0,x)\notin \mathbb{R}x^{\perp}\oplus\mathbb{R}e_3=T_xM_I$. $\Box$

\medskip

For $r>0$ let
$$
M_E(r)=\{y\in M_E:|y-e_3|<r\}.
$$

\begin{corollary}
There is a decreasing sequence $(r_j)_{j\in\mathbb{N}}$ with $r_j\to0$ as $j\to\infty$  so that for every $j\in\mathbb{N}$ there exists a continuously differentiable map
$$
t_j:M_E(r_j)\to(0,\infty)
$$
with $t_j(e_3)=t_{I,j}-t_{E,\epsilon_j}$,  $F_{\epsilon_j}(t_j(y),y)\in\{z\in M_I:P_Lz\neq-h_{\epsilon_j}(t_{I,j})\}$ for all $y\in M_E(r_j)$, and
$$
\partial_1F_{\epsilon_j}(0,y)\notin L\quad\mbox{for all}\quad y\in M_E(r_j).
$$
\end{corollary}

{\bf Proof.} Let $j\in\mathbb{N}$ be given. With the function $G:\mathbb{R}^3\to\mathbb{R}$, $G(z)=z_1^2+z_2^2-1$, the relation $F_{\epsilon_j}(t,y)\in M_I$ is equivalent to the equation $G(F_{\epsilon_j}(t,y))=0$. Using Proposition 4.2 (ii) with $x=F_{\epsilon_j}(t_{I,j}-t_{E,\epsilon_j},e_3)$  we obtain
$$
\partial_1(G\circ F_{\epsilon_j}(t_{I,j}-t_{E,\epsilon_j},e_3)=2<\partial_1F_{\epsilon_j}(t_{I,j}-t_{E,\epsilon_j},e_3),F_{\epsilon_j}(t_{I,j}-t_{E,\epsilon_j},e_3)>=2<\partial_1F_{\epsilon_j}(0,x),x)>\,\,<0.
$$
Therefore the Implicit Function Theorem applies and yields a continuously differentiable positive function $t_j^{\ast}$  on a neighbourhood $N_j$ of $e_3$ in $\mathbb{R}^3$ which satisfies $t_j^{\ast}(e_3)=t_{I,j}-t_{E,\epsilon_j}$ and $F_{\epsilon_j}(t_j^{\ast}(y),y)\in M_I$ for all $y\in N_j$. Given any real $r>0$ it follows from Proposition 4.2 (i) and continuity that there exists $r_j\in(0,r)$ with $M_E(r_j)\subset N_j$ and $\partial_1F_{\epsilon_j}(0,y)\notin L$ for all $y\in M_E(r_j)$. Using $F_{\epsilon_j}(t_j^{\ast}(e_3),e_3)=h_{\epsilon_j}(t_{I,j})\in L$ and continuity we also achieve
$P_LF_{\epsilon_j}(t^{\ast}_j(y),y)\neq-h_{\epsilon_j}(t_{I,j}))$ for all $y\in M_E(r_j)$.
Let $t_j$ be the restriction of $t_j^{\ast}$ to $M_E(r_j)$ . 

The desired decreasing sequence $(r_j)_{j\in\mathbb{N}}$ can be obtained recursively. $\Box$

\medskip

For every $j\in\mathbb{N}$ the {\it exterior map}
$$
E_j:M_E(r_j)\to \{z\in M_I:P_Lz\neq-h_{\epsilon_j}(t_{I,j})\},\quad E_j(y)=F_{\epsilon_j}(t_j(y),y),
$$
into the open subset $\{z\in M_I:P_Lz\neq-h_{\epsilon_j}(t_{I,j})\}=\{z\in M_I:P_Lz\neq-E_j(e_3))\}$ of the manifold $M_I$ is continuously differentiable. Compare Figure 2, bottom.

\begin{corollary}
Let $j\in\mathbb{N}$ be given and $x=E_j(e_3)$. Then $T_{e_3}E_j(T_{e_3}M_E)=T_xM_I$.
\end{corollary}

{\bf Proof.} Let $v=\partial_1F_{\epsilon_j}(0,e_3)$ and $w=h'_{\epsilon_j}(t_{I,j})=\partial_1F_{\epsilon_j}(t_j(e_3),e_3)=\partial F_{\epsilon_j}(0,x)$. Then $v\notin L=T_{e_3}M_E$ and $w\notin T_xM_I$. The map $T_{e_3}E_j:T_{e_3}M_E\to T_xM_I$ is given by
$$
T_{e_3}E_j(z)=P_jD_2F_{\epsilon_j}(t_j(e_3),e_3)z
$$
with the projection $P_j:\mathbb{R}^3\to\mathbb{R}^3$ along $\mathbb{R}w$ onto $T_xM_I$. 
%%% REFERENCE ? Hirsch-Smale or successor of that book ?
The isomorphism $D_2F_{\epsilon_j}(t_j(e_3),e_3)$ sends $v$ to $w$ and maps $T_{e_3}M_E$ onto a 2-dimensional space $Q$ which is complementary to $\mathbb{R}w$. The projection $P_j$ maps the complementary space $Q$ onto the complementary space $T_xM_I$. $\Box$

\section{Inner maps}

We begin with the travel time from the strip $\{x\in M_I:0<x_3<1\}$ to the plane $M_E$.

\begin{proposition}
Assume $0<\eta<e^{\sigma}$, $e^{\sigma}+\eta<1$, and $0<\epsilon<\min\{\epsilon(\eta),\epsilon_B\}$. Let $x\in M_I$ with $0<x_3<1$ be given. Then there exists $t=\tau_{\epsilon}(x)>0$ with $0<F_{\epsilon,3}(s,x)<1$ for $0\le s<t$ and
$F_{\epsilon,3}(t,x)=1$.
\end{proposition}

{\bf Proof.} 1. We show $|F_{\epsilon,3}(s,x)|>1$ for some $s>0$. Suppose this is false. Then Proposition 2.3 yields $F_{\epsilon}(n,x)\in B_1$  for all integers $n\ge0$. Using an estimate from Proposition 2.3 and $1-\eta>e^{\sigma}$ we get 
$$
|F_{\epsilon,3}(n,x)|=|P_UF_{\epsilon}(n,x)|\ge e^{un}(1-\eta)^{n+1}|P_Ux|\ge e^{(u+\sigma)n}e^{\sigma}|x_3|>0
$$
for all integers $n\ge0$, and the hypothesis $u+\sigma>0$  yields a contradiction to the assumption that $|F_{\epsilon,3}(s,x)|$, $s\ge0$, is bounded.

2. It follows that
$$
0<\inf\{s\ge0:|F_{\epsilon,3}(s,x)|\ge1\}<\infty.
$$
Set $t=\tau_{\epsilon}(x)=\inf\{s\ge0:|F_{\epsilon,3}(s,x)|\ge1\}$. Then $|F_{\epsilon,3}(s,x)|<1$ for $0\le s<t$, and
by continuity, $|F_{\epsilon,3}(t,x)|=1$. Let $n=n_{\epsilon}(x)$ denote the largest integer in $[0,t]$. Proposition 2.3 yields $F_{\epsilon}(j,x)\in B_1$ for $j=0,\ldots,n$. By $0<x_3$, $x\in B_1\setminus L$. Using Proposition 2.1 (ii) we infer $F_{\epsilon}(s,x)\in B\setminus L$ on $[0,t]$, hence $|F_{\epsilon,3}(s,x)|>0$ on $[0,t]$. Finally, $0<x_3$ and continuity combined yield $0<F_{\epsilon,3}(s,x)$ on $[0,t]$. $\Box$

\medskip

In order to use the travel time $\tau_{\epsilon}(x)$ of Proposition 5.1 for an {\it inner map} with values $F_{\epsilon}(\tau_{\epsilon}(x),x)$ in the domain $M_E(r_j)$ of an exterior map $E_j$ we observe that due to hyperbolic behavior of flowlines close to the stationary point $0$ the relation $F_{\epsilon}(\tau_{\epsilon}(x),x)\in M_E(r_j)$ should hold for $x\in M_I$ with $0<x_3$ sufficiently small. Recall $\epsilon_j$ from Proposition 4.1 and $r_j$ from Corollary 4.3. 

\begin{proposition}
Assume $0<\eta<e^{\sigma}$ and $1+\eta e^{-\sigma}<e^{-\sigma/2}$. Consider an integer $j$ so large that
$$
\epsilon_j<\min\{\epsilon(\eta),\epsilon_B\}.
$$
For every $x\in M_I$ with $0<x_3<1$ the largest integer $n=n_j(x)$ in $[0,\tau_{\epsilon_j}(x))$ satisfies
\begin{equation}
n>\frac{1}{u+\log(2)}\log\left(\frac{1}{x_3}\right)-1.
\end{equation}
For $\delta_j^{\ast}\in(0,1)$ so small that
\begin{equation}
\left(1+\frac{2}{\sigma}\log(r_j)\right)<\frac{1}{u+\log(2)}\log\left(\frac{1}{\delta_j^{\ast}}\right)-1.
\end{equation}
we have  
$$
|F_{\epsilon_j}(\tau_{\epsilon_j}(x),x)-e_3|=|P_LF_{\epsilon_j}(\tau_{\epsilon_j}(x),x)|<r_j\quad\mbox{for all}\quad x\in M_I\quad\mbox{with}\quad0<x_3<\delta_j^{\ast}.
$$
\end{proposition}

{\bf Proof.} 1. On (5). Notice that we have $e^{\sigma}+\eta<1$, so that the hypothesis concerning $\eta>0$ in Proposition 5.1 is satisfied. Let $x\in M_I$ with $0<x_3<1$ be given and let $t=\tau_{\epsilon_j}(x)$. We derive an estimate of the largest integer $n=n_j(x)$  in $[0,t)$. Proposition 5.1 in combination with Proposition 2.3, both for $\epsilon=\epsilon_j$, yield
\begin{eqnarray}
1 & = & |P_UF_{\epsilon_j}(t,x)|\le e^{u t}(1+\eta)^{n+1}|P_Ux|=e^{u t}(1+\eta)^{n+1}x_3\nonumber\\
& < & (2 e^u)^{n+1}x_3\quad\mbox{(with}\quad t\le n+1\quad\mbox{and}\quad\eta<1),\nonumber
\end{eqnarray}
hence
$$
n>\frac{1}{u+\log(2)}\log\left(\frac{1}{x_3}\right)-1.
$$

2. Now consider $\delta_j^{\ast}\in(0,1)$ which satisfies the inequality (5), and let $x\in M_I$ with $0<x_3<\delta_j^{\ast}$ be given.  Using $|P_Lx|=1$ from $x\in M_I$ and an upper estimate from Proposition 2.3 we have, with $n=n_j(x)$,
$$
|P_LF_{\epsilon_j}(t,x)|\le e^{\sigma t}(1+\eta e^{-\sigma})^{n+1}|P_Lx|\le e^{\sigma\cdot n}e^{(-\sigma/2)(n+1)}=e^{(\sigma/2)(n-1)}.
$$
Thereby the desired inequality $|P_LF_{\epsilon_j}(t,x)|<r_j$ follows from $e^{(\sigma/2)(n-1)}<r_j$
which is equivalent to
$$
n>1+\frac{2}{\sigma}\log(r_j).
$$
The preceding equation follows from the lower estimate (5) of $n=n_j(x)$ in Part 1 in combination with $x_3<\delta_j^{\ast}$ and with the smallness assumption (6) on $\delta_j^{\ast}$. $\Box$

\medskip

Next we use transversality of the flow at points of $M_E(r_j)$ as prepared in Corollary 4.3  in order to obtain smoothness of the travel time.

\begin{corollary}
Let $\eta$, $j$, and $\delta_j^{\ast}$ satisfy the hypotheses of Proposition 5.2. Then the map
$$
\tau^{\ast}_j:\{x\in M_I:0<x_3<\delta_j^{\ast}\}\to(0,\infty),\quad\tau^{\ast}_j(x)=\tau_{\epsilon_j}(x),
$$
is continuously differentiable.
\end{corollary}

{\bf Proof.} 1. Continuity.  Let $x\in M_I$ with $0<x_3<\delta_j^{\ast}$  be given, and let $t=\tau_j^{\ast}(x)=\tau_{\epsilon_j}(x)$. Let $y=F_{\epsilon_j}(t,x)$ and observe that because of $F_{\epsilon_j,3}(s,x)<1=F_{\epsilon_j,3}(t,x)$ for $0\le s<t$ we have $\partial_1F_{\epsilon_j,3}(0,y)=\partial_1F_{\epsilon_j,3}(t,x)\ge0$. Recall $y\in M_E(r_j)$ from  Proposition 5.2. The relation $\partial_1F_{\epsilon_j}(0,y)\notin L$ from Corollary 4.3 yields $\partial_1F_{\epsilon_j,3}(0,y)\neq0$. We conclude that $\partial_1F_{\epsilon_j,3}(t,x)=\partial_1F_{\epsilon_j,3}(0,y)>0$.  Now let $\rho\in(0,t)$ be given. Then for some $s\in(t,t+\rho)$, $F_{\epsilon_j,3}(s,x)>1$. By continuity there is a neighbourhood $N_1$ of $x$ in $\mathbb{R}^3$ with $F_{\epsilon_j,3}(s,z)>1$ for all $z\in N_1$. By continuity and compactness we also find a neighbourhood $N\subset N_1$ of $x$ in $\mathbb{R}^3$ so that for all $z\in N$  and for all $s\in[0,t-\rho]$ we have  $0<F_{\epsilon_j,3}(s,z)<1$. It follows that for all $z\in N\cap\{\xi\in M_I:0<\xi_3<\delta_j^{\ast}\}$ we have
$t-\rho<\tau_{\epsilon_j}(z)<s<t+\rho$,  which yields continuity of the map $t_j^{\ast}$  at $x$.

\medskip

2. We show that locally the map $\tau_j^{\ast}$ is given by continuously differentiable maps. Let $x\in M_I$ with $0<x_3<\delta_j^{\ast}$ be given and let $t=\tau_j^{\ast}(x)$.
Then $F_{\epsilon_j,3}(t,x)=1$, and
$\partial_1F_{\epsilon_j,3}(t,x)>0$, see Part 1. The Implicit Function Theorem yields an open neighbourhood $N$ of $x$ in $\mathbb{R}^3$
and $\rho>0$ and a continuously differentiable map $\tau:N\to (t-\rho,t+\rho)$ 
with $F_{\epsilon_j,3}(\tau(z),z)=1$ for all $z\in N$, and on
$(t-\rho,t+\rho)\times N$,
\begin{equation*}
F_{\epsilon_j,3}(s,z)=1\quad\mbox{if and only if}\quad s=\tau(z).
\end{equation*}
By continuity according to Part 1, there is an open neighbourhood $N_1\subset N$ of $x$ in $\mathbb{R}^3$ so that for all $z\in N_1\cap \{\xi\in M_I:0<\xi_3<\delta_j^{\ast}\}$  we have $t-\rho<\tau_j^{\ast}(z)<t+\rho$.   Recall $F_{\epsilon_j,3}(\tau_j^{\ast}(z),z)=F_{\epsilon_j,3}(\tau_{\epsilon_j}(z),z)=1$ for all $z\in M_I$ with $0<z_3<\delta_j^{\ast}$. It follows that on $N_1\cap \{\xi\in M_I:0<\xi_3<\delta_j^{\ast}\}$ we have $\tau_j^{\ast}(z)=\tau(z)$.
The restriction of $\tau$ to $N_1\cap \{\xi\in M_I:0<\xi_3<\delta_j^{\ast}\}$  is a continuously differentiable function on the open subset $N_1\cap \{\xi\in M_I:0<\xi_3<\delta_j^{\ast}\}$ of the submanifold $M_I$. $\Box$

\medskip

 In the sequel we arrange for an inner map on a subset of $\{x\in M_I:0<x_3<\delta_j^{\ast}\}$ which can be estimated in the same way as its counterpart in \cite{W2}. Upon that we will be able to follow \cite{W2} in proving existence of chaotic motion.

\medskip

For $0<\eta<\frac{e^{\sigma}}{2}$ we abbreviate
$$
m(\eta,\sigma)= \max\left\{\log(1+\eta e^{-\sigma}),\log\left(\frac{1}{1-\eta e^{-\sigma}}\right),\arcsin\left(\frac{\eta}{e^{\sigma}-\eta}\right)\right\}
$$
and consider the following hypotheses.

For $\tilde{\eta}>0$ given,
\begin{eqnarray}
0<\eta<\frac{e^{\sigma}}{2} & \mbox{and} & 1+\eta e^{-\sigma}<e^{-\sigma/2}\\
& & \mbox{(observe that this yields}\quad \eta<e^{\sigma}\quad\mbox{and}\quad e^{\sigma}+\eta<e^{\sigma/2}<1)\nonumber\\
 & \mbox{and} &  m(\eta,\sigma)<\tilde{\eta}.
\end{eqnarray}
$j\in\mathbb{N}$ is so large that 
\begin{equation}
\epsilon_j<\min\{\epsilon(\eta),\epsilon_B\}. 
\end{equation}
$\delta_j^{\ast}\in(0,1)$ satisfies (6) and 
$\delta_j\in(0,\delta_j^{\ast})$ is so small that
\begin{equation}
\frac{1}{u+\log(2)}\log\left(\frac{1}{\delta_j}\right)-1>\frac{m(\eta,\sigma)}{\tilde{\eta}-m(\eta,\sigma)}.
\end{equation}

\begin{proposition}
Let $\tilde{\eta}>0$ be given. Assume that the relations (7)-(10) hold for $\eta>0$, $j\in\mathbb{N}$, and $\delta_j$. For every $x\in M_I$ with $0<x_3<\delta_j$ the largest integer $n=n_j(x)$ in $[0,\tau_j^{\ast}(x))$ satisfies
$$
\frac{n+1}{n} m(\eta,\sigma)<\tilde{\eta}.
$$
\end{proposition}

{\bf Proof.} Using the relations (5) and (10) we get
$$
n>\frac{1}{u+\log(2)}\log\left(\frac{1}{x_3}\right)-1>\frac{1}{u+\log(2)}\log\left(\frac{1}{\delta_j}\right)-1>\frac{m(\eta,\sigma)}{\tilde{\eta}-m(\eta,\sigma)}
$$
which yields $\frac{n+1}{n}m(\eta,\sigma)<\tilde{\eta}$. $\Box$

\medskip

We are ready for the definition of the inner map $I_j$, given $\tilde{\eta}>0$ and $\eta>0$, $j\in\mathbb{N}$, and $\delta_j$ which satisfy the relations (7)-(10). As a domain for $I_j$ we take
$$
M_{I,j}=\{x\in M_I:0<x_3<\delta_j\quad\mbox{and}\quad  P_Lx\neq-E_j(e_3)\}.
$$
Here the line  given by $P_Lx= -E_j(e_3)=-h_{\epsilon_j}(t_{I,j})$ is excluded for later use,
in order to have a global parametrization of $M_{I,j}$ available. With $\tau_j(x)=\tau_j^{\ast}(x)$ on $M_{I,j}$ we set
$$
I_j(x)=F_{\epsilon_j}(\tau_j(x),x)
$$
and obtain a continuously differentiable map $I_j:M_{I,j}\to M_E$, with values in $M_E(r_j)$ according to Proposition 5.1.
Compare Figure 3, top.

\begin{figure}
	\includegraphics[page=1,scale=0.7]{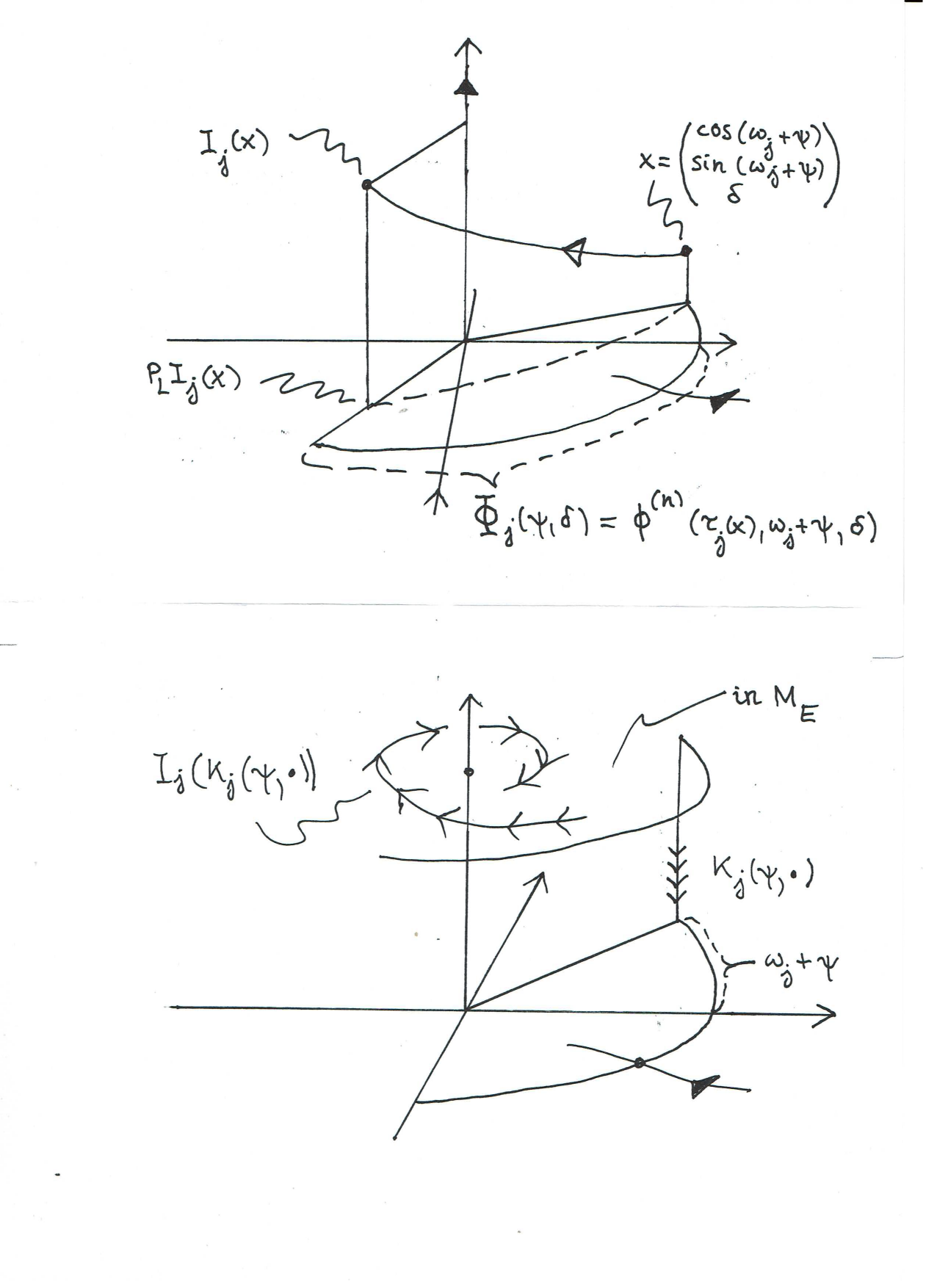}
	\caption{Top: The inner map and related angles. Bottom: The inner map along vertical line segments.}
\end{figure}

\begin{corollary}
Let $\tilde{\eta}>0$ and assume that $\eta>0$, $j\in\mathbb{N}$, and $\delta_j$ satisfy the relations (7)-(10). Let $x\in M_{I,j}$. Then we have
\begin{eqnarray*}
e^{(\sigma-\tilde{\eta})\tau_j(x)} & \le|P_LI_j(x)|\le & e^{(\sigma+\tilde{\eta})\tau_j(x)},\\
e^{(u-\tilde{\eta})\tau_j(x)}x_3 & \le 1\le & e^{(u+\tilde{\eta})\tau_j(x)}x_3.
\end{eqnarray*}
For $\psi\in\mathbb{R}$  with
$$
x=\left(\begin{array}{c}\cos(\psi)\\ \sin(\psi)\\x_3 \end{array}\right)
$$
and $\delta=x_3$ the function $\phi^{(n)}$ obtained in Proposition 3.2, with $n=n_j(x)$,  satisfies
$$
\psi-\tau_j(x)(\mu+\tilde{\eta})\le\phi^{(n)}(\tau_j(x),\psi,\delta)\le\psi-\tau_j(x)(\mu-\tilde{\eta}).
$$
\end{corollary} 

{\bf Proof.}  Proposition 5.4 shows that the hypotheses on the integer $n$ in Proposition 2.4 and Corollary 3.3 are satisfied for $n=n_j(x)$. Apply Proposition 2.4 to $t=\tau_j(x)-n_j(x)\in[0,1]$, with $|P_Lx|=1$ and $|P_Ux|=x_3$ and $P_UI_j(x)=e_3$. This yields the estimates of $P_LI_j(x)$ and $P_UI_j(x)$. 

From $|P_Lx|=1$, $x\in B_1\setminus U$.  An application of Corollary 3.3 to $t=\tau_j(x)-n_j(x)\in[0,1]$ yields the estimate of the
angle function $\phi^{(n)}$ with $n=n_j(x)$. $\Box$

\begin{corollary}
Let $\tilde{\eta}>0$ be given with $\tilde{\eta}<u$ and assume that $\eta>0$, $j\in\mathbb{N}$, and $\delta_j$ satisfy the relations (7)-(10). For every $x\in M_{I,j}$,
$$
\frac{1}{u+\tilde{\eta}}\log\left(\frac{1}{x_3}\right)\le\tau_j(x)\le\frac{1}{u-\tilde{\eta}}\log\left(\frac{1}{x_3}\right).
$$
\end{corollary}

\section{The return map in the plane}

We begin with parametrizations of the open subsets $\{x\in M_I:P_Lx\neq-E_j(e_3)\}$ of the submanifold $M_I$. Let $j\in\mathbb{N}$ be given. Consider $\omega_j\in[-\pi,\pi)$ determined by
$$
\left(\begin{array}{c}\cos(\omega_j)\\ \sin(\omega_j)\\0\end{array}\right)=E_j(e_3)\quad(=h_{\epsilon_j}(t_{I,j})).
$$ 
The map
$$
K_j:(-\pi,\pi)\times\mathbb{R}\to \{x\in M_I:P_Lx\neq-E_j(e_3)\},\quad K_j(\psi,\delta)=\left(\begin{array}{c}\cos(\omega_j+\psi)\\ \sin(\omega_j+\psi)\\ \delta\end{array}\right),
$$
is a continuously differentiable diffeomorphism with $K_j(0,0)=E_j(e_3)$. The return map $R_j=E_j\circ I_j$ sends its domain $M_{I,j}=\{x\in M_I:0<x_3<\delta_j,P_Lx\neq-E_j(e_3)\}$ into the set $\{z\in M_I:P_Lz\neq-E_j(e_3)\}$ which equals the image $K_j((-\pi,\pi)\times\mathbb{R})$. By the {\it return map in the plane} we mean the continuously differentiable map
$$
Q_j:(-\pi,\pi)\times(0,\delta_j)\to(-\pi,\pi)\times\mathbb{R}
$$
given by $Q_j(\psi,\delta)=K_j^{-1}(R_j(K_j(\psi,\delta)))$.

We also need information about a coordinate representation of the exterior map $E_j$ alone. Corollary 4.4 yields that the derivative $T_{e_3}(K_j^{-1}\circ E_j)$ is an isomorphism from $L=T_{e_3}M_E$ onto the plane $\mathbb{R}^2$. So it sends basis vectors $v_j$ and $w_j$ of $L$ to the vectors $\left(\begin{array}{c}1\\0\end{array}\right)$ and $\left(\begin{array}{c}0\\1\end{array}\right)$, respectively. 

Let $\kappa_j$ denote the isomorphism $L\to\mathbb{R}^2$ given by
$$
\kappa_j(\xi v_j+\eta w_j)=\left(\begin{array}{c}\xi\\ \eta\end{array}\right).
$$
The restriction $P_j$ of $\kappa_j\circ P_L$ to the open subset $M_E(r_j)$ of the submanifold $M_E$ defines a continuously differentiable diffeomorphism onto an open neighbourhood of $0=P_j(e_3)$ in $\mathbb{R}^2$. Obviously,
$$
P_j(e_3+\xi v_j+\eta w_j)=\left(\begin{array}{c}\xi\\ \eta\end{array}\right)\quad\mbox{for all reals}\quad\xi,\eta.
$$
As $T_{e_3}P_j$ and $T_{e_3}(K_j^{-1}\circ E_j)$ act the same on the basis $v_j,w_j$ of $T_{e_3}M_E=L$ the {\it exterior map in coordinates}
$$
K_j^{-1}\circ E_j\circ P_j^{-1}:\mathbb{R}^2\supset P_j(M_E(r_j))\to\mathbb{R}^2
$$
satisfies
\begin{equation}
	D(K_j^{-1}\circ E_j\circ P_j^{-1})(0)=id_{\mathbb{R}^2}.
\end{equation}

\begin{corollary}
Let $j\in\mathbb{N}$ and $\beta>0$ be given. There exists $\alpha_j=\alpha_j(\beta)\in(0,\pi)$ so that for all $(\psi,\delta)\in[-\alpha_j,\alpha_j]\times[-\alpha_j,\alpha_j]$ we have
\begin{eqnarray}
(\psi,\delta) & \in & P_j(M_E(r_j),\nonumber\\
|(K_j^{-1}\circ E_j\circ P_j^{-1})(\psi,\delta)-(\psi,\delta)| & \le & \beta|(\psi,\delta)|,\\
|D(K_j^{-1}\circ E_j\circ P_j^{-1})(\psi,\delta)-id_{\mathbb{R}^2}|
& \le & \beta.
\end{eqnarray}
\end{corollary}

We proceed to estimates of the range of the {\it inner map in coordinates} 
$$
(-\pi,\pi)\times(0,\delta_j)\ni(\psi,\delta)\mapsto P_j(I_j(K_j(\psi,\delta)))\in\mathbb{R}^2
$$
and of the return map in the plane $Q_j$. 

\begin{proposition}
Assume $0<\tilde{\eta}<-\sigma/2$ and consider $\eta>0, j\in\mathbb{N},\delta_j$ so that the relations (7-10) hold. 
Let $\beta\in(0,\frac{1}{2}]$ be given. Consider $\alpha_j=\alpha_j(\beta)>0$ according to Corollary 6.1, with $\alpha_j<\delta_j$. Then
$$
\delta_{\beta,j}=\left(\frac{2}{3(|\kappa_j|+1)}\alpha_j\right)^{\frac{3u}{-\sigma}}
$$ 
satisfies $\delta_{\beta,j}\le\frac{2}{3}\alpha_j$,
and for all $(\psi,\delta)\in[-\alpha_j,\alpha_j]\times(0,\delta_{\beta,j}]$ we have
\begin{eqnarray}
|P_j(I_j(K_j(\psi,\delta)))| & \le \frac{2}{3}\alpha_j,\\
Q_j(\psi,\delta) & \in & [-\alpha_j,\alpha_j]\times[-\alpha_j,\alpha_j].
\end{eqnarray}
\end{proposition}

{\bf Proof.} 1. In order to show $\delta_{\beta,j}\le\frac{2}{3}\alpha_j$ notice first that
$$
\delta_{\beta,j}<\alpha_j<\delta_j<1.
$$
The hypothesis on $\tilde{\eta}$ yields
$$
\frac{u+\tilde{\eta}}{-\sigma-\tilde{\eta}}<\frac{\frac{3u}{2}}{\frac{-\sigma}{2}}=\frac{3u}{-\sigma}
$$
It follows that
$$
\delta_{\beta,j}=\left(\frac{2}{3(|\kappa_j|+1)}\alpha_j\right)^{\frac{3u}{-\sigma}}\le\left(\frac{2}{3(|\kappa_j|+1)}\alpha_j\right)^\frac{u+\tilde{\eta}}{-\sigma-\tilde{\eta}},
$$
hence
$$
\delta_{\beta,j}^{\frac{-\sigma-\tilde{\eta}}{u+\tilde{\eta}}}\le\frac{2}{3(|\kappa_j|+1)}\alpha_j
$$
Consequently, with $0<-\sigma-\tilde{\eta}<u+\tilde{\eta}$ and $\delta_{\beta,j}<1$,
\begin{eqnarray}
\delta_{\beta,j} & \le & \delta_{\beta,j}^\frac{-\sigma-\tilde{\eta}}{u+\tilde{\eta}}\le\frac{2}{3(|\kappa_j|+1)}\alpha_j\\
& \le &\frac{2}{3}\alpha_j.\nonumber
\end{eqnarray}

2. From Corollary 6.1  with $0<\beta\le1/2$ we get $\alpha_j<\delta_j$ 
so that for all $(\tilde{\psi},\tilde{\delta})\in[-\alpha_j,\alpha_j]\times[-\alpha_j,\alpha_j]$, 
we have	
$$
|(K_j^{-1}\circ E_j\circ P_j^{-1})(\tilde{\psi},\tilde{\delta})|\le(1+\beta)|(\tilde{\psi},\tilde{\delta})|\le\frac{3}{2}|(\tilde{\psi},\tilde{\delta})|.
$$

3. Proof of (14). For  $(\psi,\delta)\in[-\alpha_j,\alpha_j]\times(0,\delta_{\beta,j}]$
we have $K_{j,3}(\psi,\delta)=\delta\in(0,\delta_{\beta,j}]\subset(0,\delta_j)$. With $x=K_j(\psi,\delta)$, Corollary 5.5 yields
$$
|P_j(I_j(K_j(\psi,\delta)))|=|\kappa_jP_LI_j(x)|\le|\kappa_j|e^{(\sigma+\tilde{\eta})\tau_j(x)}.
$$
Notice that $\tilde{\eta}<u$. Using the lower estimate of $\tau_j(x)$ from Corollary 5.6 we infer
\begin{eqnarray}
|P_j(I_j(K_j(\psi,\delta)))| & \le & |\kappa_j|\left(\frac{1}{x_3}\right)^{\frac{\sigma+\tilde{\eta}}{u+\tilde{\eta}}}\nonumber\\
& \le & |\kappa_j|\delta_{\beta,j}^{\frac{-\sigma-\tilde{\eta}}{u+\tilde{\eta}}}\le\frac{2}{3}\alpha_j\quad\mbox{(with (16))}.\nonumber
\end{eqnarray}

4. Proof of (15). For $(\psi,\delta)$ as  in Part 3 let $(\tilde{\psi},\tilde{\delta})=(P_j(I_j(K_j)))(\psi,\delta)$. Then
$$
|(\tilde{\psi},\tilde{\delta})|\le\frac{2}{3}\alpha_j,
$$
hence $(\tilde{\psi},\tilde{\delta})\in[-\alpha_j,\alpha_j]\times[-\alpha_j,\alpha_j]$, which according to Part 2 yields 
$$
|(K^{-1}_j\circ E_j\circ P_j^{-1})(\tilde{\psi},\tilde{\delta})|\le \frac{3}{2}|(\tilde{\psi},\tilde{\delta})|.
$$
It follows that
\begin{eqnarray}
|Q_j(\psi,\delta)| & = & 
|K^{-1}_j(R_j(K_j(\psi,\delta)))|=|(K^{-1}_j\circ E_j\circ P_j^{-1}\circ P_j\circ I_j)(K_j(\psi,\delta))|\nonumber\\
& = & |(K^{-1}_j\circ E_j\circ P_j^{-1})(\tilde{\psi},\tilde{\delta})|\nonumber\\
& \le & \frac{3}{2}|(\tilde{\psi},\tilde{\delta})|\le\alpha_j.\nonumber
\end{eqnarray}
Finally, use that the disk of radius $\alpha_j$ and center $0\in\mathbb{R}^2$ is contained in the square $[-\alpha_j,\alpha_j]\times[-\alpha_j,\alpha_j]$. $\Box$ 

\medskip

The last result of this section concerns continuity of the angle corresponding to $P_L I_j(K_j(\psi,\delta))$, as a function of $(\psi,\delta)\in (-\pi,\pi)\times(0,\delta_j)$.
 
\begin{proposition}
Assume $0<\tilde{\eta}<-\sigma/2$ and consider $\eta>0, j\in\mathbb{N},\delta_j$ so that the relations (7-10) hold. 
Then the function $\Phi_j:(-\pi,\pi)\times(0,\delta_j)\to\mathbb{R}$ given by
$$
\Phi_j(\psi,\delta)=\phi^{(n)}(\tau_j(x),\omega_j+\psi,\delta)
$$
with $\phi^{(n)}$ according to Proposition 3.2,  $n=n_j(x)$ the largest integer in $[0,\tau_j(x))$ and  $x=K_j(\psi,\delta)$,
satisfies 
\begin{equation}
\frac{1}{|P_LI_j(K_j(\psi,\delta))|}P_LI_j(K_j(\psi,\delta))=\left(\begin{array}{c}\cos(\Phi_j(\psi,\delta))\\ \sin(\Phi_j(\psi,\delta))\\0\end{array}\right)\quad\mbox{for all}\quad(\psi,\delta)\in(-\pi,\pi)\times(0,\delta_j)
\end{equation}
and is continuous.
\end{proposition}

{\bf Proof.} 1. The definition of $\Phi_j$ makes sense because for $(\psi,\delta)\in(-\pi,\pi)\times(0,\delta_j$ and $x=K_j(\psi,\delta)$ we have $|P_UF_{\epsilon_j}(t,x)|<1$ on $[0,\tau_j(x))$, hence $|P_UF_{\epsilon_j}(\nu,x)|<1$ for $\nu=0,\ldots,n_j(x)$, which in combination with $n_j(x)<\tau_j(x)\le n_j(x)+1$ yields
$(\tau_j(x),\omega_j+\psi,\delta)\in dom_{n_j(x)}$.

Eq. (17) holds because for $\phi=\phi^{n_j(x)}(\tau_j(x),\omega_j+\psi,\delta)=\Phi_j(\psi,\delta)$ we have
$$
\frac{1}{|P_LI_j(K_j(\psi,\delta))|}P_LI_j(x)=\frac{1}{|P_LF_{\epsilon_j}(\tau_j(x),x)|}P_LF_{\epsilon_j}(\tau_j(x),x)=\left(\begin{array}{c}\cos(\phi)\\ \sin(\phi)\\0\end{array}\right)
$$ 
due to Proposition 3.2 (i).

\medskip

2. Proof that $\Phi_j$ is continuous at  $(\psi,\delta)\in(-\pi,\pi)\times(0,\delta_j)$. Let $x=K_j(\psi,\delta)$. We have $n_j(x)<\tau_j(x)\le n_j(x)+1<n_j(x)+2$. Let a sequence $(\psi_m,\delta_ {(m)})_{m\in\mathbb{N}}$ in $(-\pi,\pi)\times(0,\delta_j)$ be given which converges to $(\psi,\delta)$. It is enough to find a subsequence so that $\Phi_j(\psi_{m_{\mu}},\delta_{(m_{\mu})})\to\Phi_j(\psi,\delta)$ as $\mu\to\infty$. 

2.1. In case $\tau_j(x)<n+1$ the continuity of $\tau_j$ 
yields an integer $m_x\ge0$ so that for all indices $m\ge m_x$ we have $n_j(x)<\tau_j(x_m)<n_j(x)+1$ for $x_m=K_j(\psi_m,\delta_{(m)})$. Hence $n_j(x)=n_j(x_m)$ for $m\ge m_x$. Consequently, $(\tau_j(x_m),\omega_j+\psi_m,\delta_{(m)})\in dom_{n_j(x_m)}=dom_{n_j(x)}$ and
$$
\Phi_j(\psi_m,\delta_{(m)})= \phi^{(n_j(x_m))}(\tau_j(x_m),\omega_j+\psi_m,\delta_{(m)})=\phi^{(n_j(x))}(\tau_j(x_m),\omega_j+\psi_m,\delta_{(m)})\quad\mbox{for}\quad m\ge m_x.
$$
As $\phi^{(n_j(x))}$ is continuous according to Proposition 3.2 we arrive at 
$$
\lim_{m\to\infty}\Phi_j(\psi_m,\delta_{(m)})=\lim_{m\to\infty}\phi^{(n_j(x))}(\tau_j(x_m),\omega_j+\psi_m,\delta_{(m)})=\phi^{(n_j(x))}(\tau_j(x),\omega_j+\psi,\delta)=\Phi_j(\psi,\delta).
$$

2.2. In case $\tau_j(x)=n_j(x)+1$ we have $|P_UF_{\epsilon_j}(n_j(x)+1,x)|=|P_UF_{\epsilon_j}(\tau_j(x),x)|=1$ in addition to $|P_UF_{\epsilon_j}(\nu,x)|<1$ for $\nu=0,\ldots,n_j(x)$ and conclude that $(\tau_j(x),\omega_j+\psi,\delta)\in dom_{n_j(x)+1}$. 

By Proposition 3.2 (ii), $\phi^{(n_j(x)+1)}(n_j(x)+1,\omega_j+\psi,\delta)=\phi^{(n_j(x))}(n_j(x)+1,\omega_j+\psi,\delta)$.

We distinguish the subcases that the indices $m$ with $\tau_j(x_m)\le n_j(x)+1$ are bounded or not. 

2.2.1. If the indices with 
$\tau_j(x_m)\le n_j(x)+1$ are unbounded then there is a strictly increasing sequence $(m_{\mu})_{\mu\in\mathbb{N}}$ of positive integers with $\tau_j(x_{m_{\mu}})\le n_j(x)+1$ for all $\mu\in\mathbb{N}$. As in Part 2.1 we find
$\Phi_j(\psi_{m_{\mu}},\delta_{(m_{\mu})})\to \Phi_j(\psi,\delta)$
for $\mu\to\infty$.

2.2.2. If there is an upper bound $m_x\in\mathbb{N}$ for the indices $m$ with $\tau_j(x_m)\le n_j(x)+1$ then $n_j(x)+1<\tau_j(x_m)$ for all indices $m>m_x$. In addition we may assume $\tau_j(x_m)<n_j(x)+2$ for all $m>m_x$. It follows that  $n_j(x)+1=n_j(x_m)$ and
$$
\Phi_j(\psi_m,\delta_{(m)})= \phi^{(n_j(x_m))}(\tau_j(x_m),\omega_j+\psi_m,\delta_{(m)})=\phi^{(n_j(x)+1)}(\tau_j(x_m),\omega_j+\psi_m,\delta_{(m)})
$$
for $m>m_x$. Using continuity of $\phi^{(n_j(x)+1)}$ from Proposition 3.2 we find
\begin{eqnarray*}
\lim_{m\to\infty}\Phi_j(\psi_m,\delta_{(m)}) & = & \lim_{m\to\infty}\phi^{(n_j(x)+1)}(\tau_j(x_m),\omega_j+\psi_m,\delta_{(m)})\\
& = & \phi^{(n_j(x)+1)}(\tau_j(x),\omega_j+\psi,\delta)=\phi^{(n_j(x)+1)}(n_j(x)+1,\omega_j+\psi,\delta)\\
& = & \phi^{(n_j(x))}(n_j(x)+1,\omega_j+\psi,\delta)
=\phi^{(n_j(x))}(\tau_j(x),\omega_j+\psi,\delta)=\Phi_j(\psi,\delta).
\quad\Box
\end{eqnarray*}

\section{Curves expanded by the return map in the plane}

In this section we consider curves $g:[a,b]\to(-\pi,\pi)\times(0,\delta_j)$ which connect level sets $(-\pi,\pi)\times\{\Delta_1\}$ and $(-\pi,\pi)\times\{\Delta_2\}$ with $0<\Delta_1<\Delta_2<\delta_j$. We find subintervals $[a_0,b_0]$ and $[a_1,b_1]$ so that the angle function $\Phi_j$ sends $g((a_0,b_0))$ and $g((a_1,b_1))$ into disjoint sets $M_0$ and $M_1$, and transport by the return map in the plane $Q_j$ yields
two curves which again connect the said level sets. Be aware that the angle function $\Phi_j$ is not directly related to the return map in the plane but only to the inner map in coordinates (depicted in Figure 3, bottom), which is the first composite of the return map in the plane.

\medskip

Throughout this section we assume $0<\tilde{\eta}<\min\{\mu,-\sigma/2\}$, and that $\eta>0$, $j\in\mathbb{N}$, and $\delta_j\in(0,1)$ satisfy the relations (7-10). We begin with a comparison of angles $\Phi_j(\psi,\delta)$ for arguments in level sets as above. Let
\begin{eqnarray*}
	c & = & c_{\tilde{\eta}}=\frac{(u+\tilde{\eta})(\mu+\tilde{\eta})}{(u-\tilde{\eta})(\mu-\tilde{\eta})}\\
	& > & 1,\\
	k & = & k_{\tilde{\eta}}=e^{-6\pi\frac{u+\tilde{\eta}}{\mu-\tilde{\eta}}}\\
	& < & 1.
\end{eqnarray*}
For $\Delta_2\in(0,\delta_j)$ we set
\begin{eqnarray*}
\Delta_1 & = & \Delta_1(\Delta_2)=k\Delta_2^{c}\\
& < & \Delta_2.
\end{eqnarray*}

\begin{proposition}
Let $\Delta_2\in(0,\delta_j)$ be given and consider $\Delta_1=\Delta_1(\Delta_2)=k\Delta_2^{c}$. Then
$$
4\,\pi\le\Phi_j(\psi,\Delta_2)-\Phi_j(\gamma,\Delta_1)\quad\mbox{for all}\quad\psi,\gamma\quad\mbox{in}\quad(-\pi,\pi).
$$
\end{proposition}

{\bf Proof.} Assume $-\pi<\psi<\pi$, $-\pi<\gamma<\pi$. Recall the definition of $\Phi_j$ in Proposition 6.3 and apply the last estimate in Corollary 5.5, and the estimate of the travel time in Corollary 5.6. This yields the inequalities
$$
\Phi_j(\psi,\Delta_2)-(\omega_j+\psi)\ge
(-\mu-\tilde{\eta})\cdot\frac{1}{u-\tilde{\eta}}\log\left(\frac{1}{\Delta_2}\right)
$$
and
$$
\Phi_j(\gamma,\Delta_1)-(\omega_j+\gamma)\le(-\mu+\tilde{\eta})\cdot\frac{1}{u+\tilde{\eta}}\log\left(\frac{1}{\Delta_1}
\right)
$$
from which we obtain
\begin{eqnarray*}
\Phi_j(\psi,\Delta_2)-\Phi_j(\gamma,\Delta_1) & \ge & -2\pi+\frac{\mu+\tilde{\eta}}{u-\tilde{\eta}}\log(\Delta_2)
	-\frac{\mu-\tilde{\eta}}{u+\tilde{\eta}}\log(\Delta_1)\\
& 	\ge &-2\pi+\frac{\mu+\tilde{\eta}}{u-\tilde{\eta}}\log(\Delta_2)
-\frac{\mu-\tilde{\eta}}{u+\tilde{\eta}}[\log(k)+c\log(\Delta_2)]\\
& = & -2\pi+6\pi+\left(\frac{\mu+\tilde{\eta}}{u-\tilde{\eta}}-c\frac{\mu-\tilde{\eta}}{u+\tilde{\eta}}\right)\log(\Delta_2)=4\pi.\quad\Box
\end{eqnarray*}

From here on let also $\beta\in(0,1/2]$ be given and consider $\alpha_j=\alpha_j(\beta)\in(0,\delta_j)$ according to Proposition 6.2.  Proposition 7.1 shows that the quantities
$$
m_1=m_1(j,\Delta_1)=\max_{|\gamma|\le\alpha_j}\Phi_j(\gamma,\Delta_1)\quad\mbox{and}\quad m_2=m_2(j,\Delta_2)=\min_{|\psi|\le\alpha_j}\Phi_j(\psi,\Delta_2)
$$ 
satisfy $m_1+4\pi\le m_2$. Also, there exists $\psi_j\in[m_1+\pi,m_2-\pi]$ with
$$
\left(\begin{array}{c}\cos(\psi_j)\\ \sin(\psi_j)\\0\end{array}\right)=\frac{1}{|w_j|}w_j.
$$

\begin{proposition}
(Angles along curves connecting vertical levels). Consider $\Delta_2\in(0,\delta_j)$ and $\Delta_1=\Delta_1(\Delta_2)$ as above. Let a curve $g:[a,b]\to(-\pi,\pi)\times(0,\delta_j)$ be given with $g_2(b)=\Delta_2$ and $g_2(a)=\Delta_1$. Then there exist  
$a'_0<b'_0\le a'_1<b'_1$ in $[a,b]$
such that
$$
\Phi_j(g(t))\in(\psi_j-\pi,\psi_j)\quad\mbox{on}\quad(a'_0,b'_0),\quad\Phi_j(g(a'_0))=
\psi_j-\pi,\quad\Phi_j(g(b'_0))=\psi_j,
$$
$$
\Phi_j(g(t))\in(\psi_j,\psi_j+\pi)\quad\mbox{on}\quad(a'_1,b'_1),\quad\Phi_j(g(a'_1))=\psi_j,\quad\Phi_j(g(b'_1))=\psi_j+\pi.
$$
\end{proposition}

Compare Figure 4.

\medskip

\begin{figure}
	\includegraphics[page=1,scale=0.7]{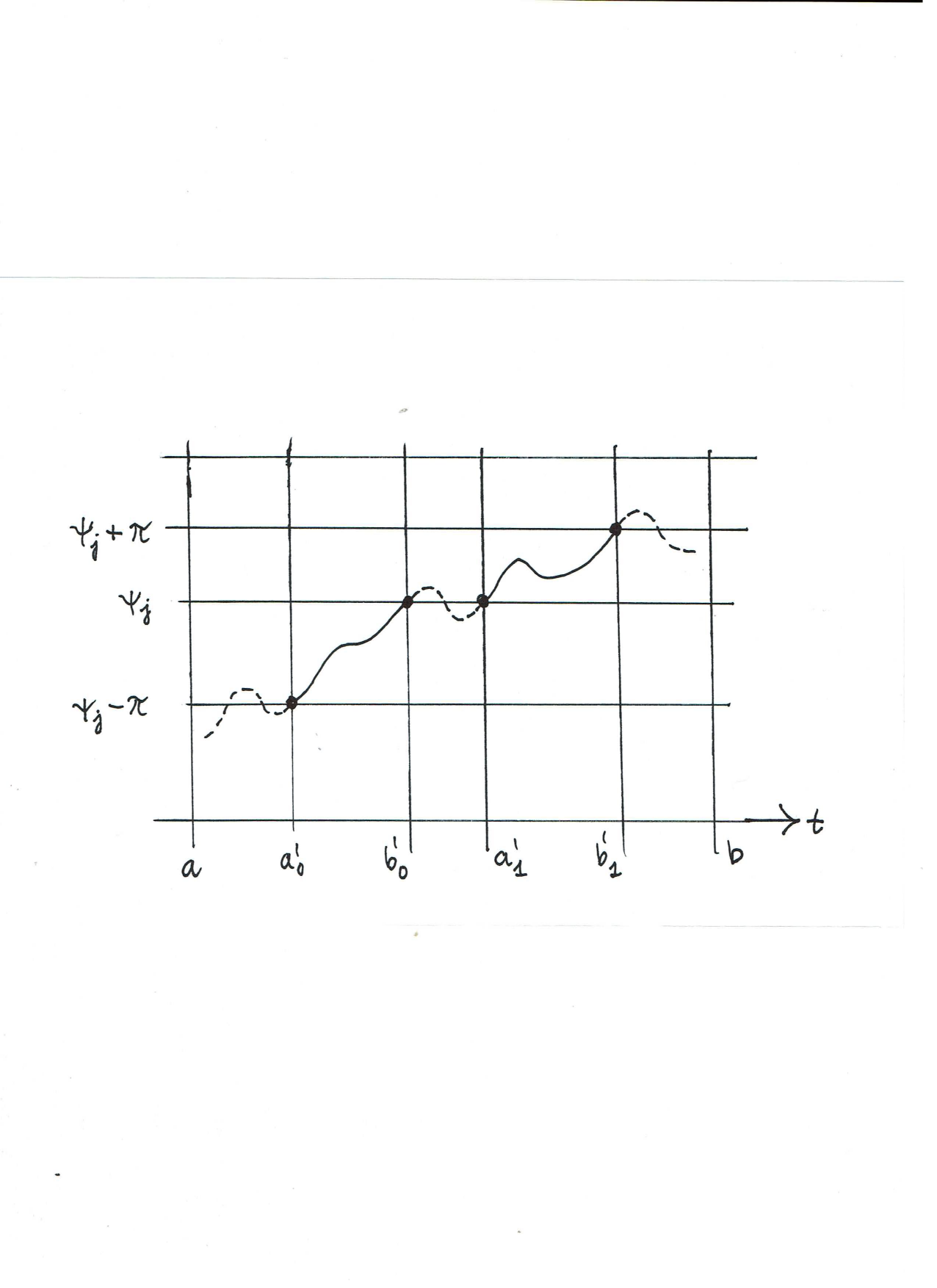}
	\caption{The map $\Phi_j\circ g$}
\end{figure}

{\bf Proof.} 1. We construct $a'_1$ and $b'_1$. From
$$
\Phi_j(g(a))\le m_1<m_1+\pi\le\psi_j\le m_2-\pi<m_2\le\Phi_j(g(b))
$$
we have
$$
\Phi_j(g(a))\le\psi_j-\pi<\psi_j<\psi_j+\pi\le\Phi_j(g(b)).
$$
By continuity, $\psi_j=\Phi_j(g(t))$ for some $t\in(a,b)$. Again by continuity there exists $b'_1\in(t,b]$ with 
$\Phi_j(g(s))<\psi_j+\pi$ on $[t,b'_1)$ and 
$\Phi_j(g(b'_1))=\psi_j+\pi$. Upon that, there exists
$a'_1\in[t,b'_1)$ with $\psi_j<\Phi_j(g(s))$ on $(a'_1,b'_1]$ and $\Phi_j(g(a'_1))=\psi_j$.

\medskip
2. The construction of $a'_0$ and $b'_0$ with $b'_0\le a'_1$ is analogous. $\Box$

\medskip

We turn to the position of $Q_j(\psi,\delta)$ for arguments $(\psi,\delta)\in(-\pi,\pi)\times[\Delta_1,\Delta_2]$. A look at Eq. (17) in Proposition 6.3
confirms that in the cases
$$
\Phi_j(\psi,\delta)=\psi_j-\pi,\quad\Phi_j(\psi,\delta)=\psi_j,\quad\Phi_j(\psi,\delta)=\psi_j+\pi
$$
the value $P_LI_j(K_j(\psi,\delta))$ belongs to the rays
$$
(0,\infty)(-w_j),\quad(0,\infty)w_j,\quad(0,\infty)(-w_j),\quad\mbox{respectively,}
$$ 
hence $P_j(I_j(K_j(\psi,\delta)))$ is on the vertical axis in $\mathbb{R}^2$.

\begin{proposition}
(From angles to vertical levels) Assume  in addition to the hypotheses made in the present section up to here that $\tilde{\eta}$ satisfies
\begin{equation}
c_{\tilde{\eta}}\frac{-\sigma+\tilde{\eta}}{u-\tilde{\eta}}<1.
\end{equation} 
Let $\beta\in(0,1/2]$ be given and choose reals $\alpha_j=\alpha_j(\beta)\in(0,\delta_j)$ and $\delta_{\beta,j}\in
(0,2\alpha_j/3]$ according to Proposition 6.2.
Consider $\Delta_2\in(0,\delta_{\beta,j})$ so small that
\begin{equation}
2\sqrt{2}\Delta_2<\frac{1}{|\kappa_j^{-1}|}k^{\frac{-\sigma+\tilde{\eta}}{u-\tilde{\eta}}}\Delta_2^{c\frac{-\sigma+\tilde{\eta}}{u-\tilde{\eta}}},\quad\mbox{with}\quad c=c_{\tilde{\eta}}\quad\mbox{and}\quad k=k_{\tilde{\eta}}.
\end{equation}
Let $\Delta_1=\Delta_1(\Delta_2)$, $(\psi,\delta)\in[-\alpha_j,\alpha_j]\times[\Delta_1,\Delta_2]$, and $z=Q_j(\psi,\delta)$. Then

\medskip

(i)  $|z|>\sqrt{2}\Delta_2$.

\medskip

(ii) In the cases $\Phi_j(\psi,\delta)=\psi_j-\pi$, $\quad\Phi_j(\psi,\delta)=\psi_j$, $\quad\Phi_j(\psi,\delta)=\psi_j+\pi$,
we have 
$$
z_2<-\Delta_2,\quad z_2>\Delta_2,\quad
z_2<-\Delta_2,\quad\mbox{respectively}.
$$
\end{proposition}

Compare Figure 5.

\medskip

\begin{figure}
	\includegraphics[page=1,scale=0.7]{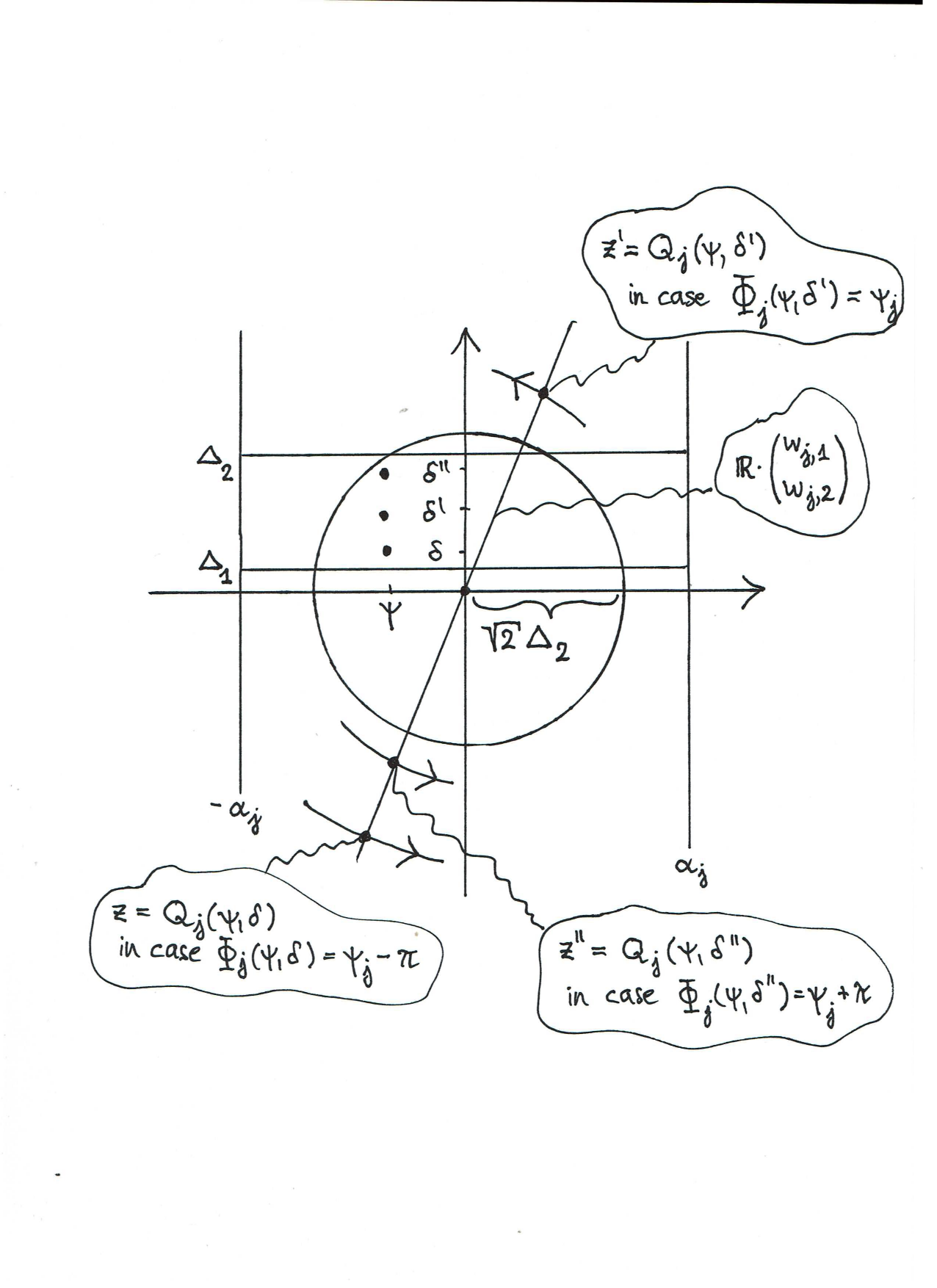}
	\caption{Positions of $Q_j(\psi,\delta)$ depending on the angle $\Phi_j(\psi,\delta)$}
\end{figure}

{\bf Proof.} 1. On assertion (i). Observe first that due to Proposition 6.2  the rectangle 
$[-\alpha_j,\alpha_j]\times[\Delta_1,\Delta_2]$ 
is contained in the domain of definition of the maps $P_j(I_j(K_j(\cdot,\cdot)))$ and $Q_j$.
Let $(\psi,\delta)\in[-\alpha_j,\alpha_j]\times[\Delta_1,\Delta_2]$ be given, and let 
$$
x=P_j(I_j(K_j(\psi,\delta)))\quad(=\kappa_jP_LI_j(K_j(\psi,\delta))).
$$
Observe $K_{j,3}(\psi,\delta)=\delta$ and apply Corollaries 5.5 and 5.6. This yields
$$
|x|\ge\frac{1}{|\kappa_j^{-1}|}\delta^{\frac{-\sigma+\tilde{\eta}}{u-\tilde{\eta}}}\ge
\frac{1}{|\kappa_j^{-1}|}\Delta_1^{\frac{-\sigma+\tilde{\eta}}{u-\tilde{\eta}}}
=\frac{1}{|\kappa_j^{-1}|}k^{\frac{-\rho+\tilde{\eta}}{u-\tilde{\eta}}}\Delta_2^{c\frac{-\sigma+\tilde{\eta}}{u-\tilde{\eta}}}>2\sqrt{2}\Delta_2.
$$
We have
$$
z=Q_j(\psi,\delta)=(K_j^{-1}\circ E_j\circ I_j)(K_j(\psi,\delta))=(K_j^{-1}\circ E_j\circ P_j^{-1}\circ P_j\circ I_j)(K_j(\psi,\delta)) =(K_j^{-1}\circ E_j\circ P_j^{-1})(x).
$$
The inequality (14) in Proposition 6.2 gives us $|x|\le\frac{2}{3}\alpha_j$, hence $x\in[-\alpha_j,\alpha_j]\times[-\alpha_j,\alpha_j]$. Using the relation (12) in Corollary 6.1 we infer $|z-x|\le\beta|x|$. It follows that
$$
|z|\ge|x|-\beta|x|\ge\frac{1}{2}|x|>\sqrt{2}\Delta_2.
$$

2. On assertion (ii) for the case $(\psi,\delta)\in[-\alpha_j,\alpha_j]\times[\Delta_1,\Delta_2]$ with $\Phi_j(\psi,\delta)=\psi_j$. Let $z=Q_j(\psi,\delta)\in\mathbb{R}^2$.

\medskip

2.1. From (17) in Proposition 6.3 we have that $P_LI_j(K_j(\psi,\delta))$ is a positive multiple of
$$
\left(\begin{array}{c}\cos(\Phi_j(\psi,\delta))\\ \sin(\Phi_j(\psi,\delta))\\0\end{array}\right)=\left(\begin{array}{c}\cos(\psi_j)\\ \sin(\psi_j)\\0\end{array}\right)\in(0,\infty)w_j,
$$
hence 
$$
P_j(I_j(K_j(\psi,\delta)))=\kappa_j P_LI_j(K_j(\psi,\delta))\in(0,\infty)\left(\begin{array}{c}0\\1\end{array}\right).
$$ 

2.2. From Part 1, $|z|>\sqrt{2}\Delta_2$. For $x=P_j(I_j(K_j(\psi,\delta)))=\kappa_jP_LI_j(K_j(\psi,\delta))$  Part 2.1 yields $x=\left(\begin{array}{c}0\\x_2\end{array}\right)$ with $x_2>0$.

\medskip

2.3. Proof of $|z|\le\sqrt{2}z_2$: We have $|x|=x_2$. From 
$x_2-z_2\le|x_2-z_2|\le|z-x|\le\beta|x|=\beta x_2$, $z_2\ge(1-\beta)x_2>0$. Also, from $x_1=0$, $|z_1|\le|x_1|+\beta|x|=\beta x_2$. It follows that
$$
|z|^2=z_1^2+z_2^2\le\beta^2x_2^2+z_2^2\le\frac{\beta^2}{(1-\beta)^2}z_2^2+z_2^2\le 2z_2^2.
$$

2.4. Consequently, $z_2=|z_2|\ge\frac{1}{\sqrt{2}}|z|>\Delta_2$.

\medskip

3. The proofs of assertion (ii) in the two remaining cases are analogous, making use of the fact that in both cases we have that
$P_j(I_j(K_j(\psi,\delta)))=\kappa_j P_LI_j(K_j(\psi,\delta))$ is a positive multiple of 
$\left(\begin{array}{c}0\\-1\end{array}\right)$. $\Box$

\medskip

The next result makes precise what was briefly announced at the begin of the section. The disjoint sets mentioned there will be given in terms of the angle $\Phi_j(\psi,\delta)$ corresponding to the value $P_j(I_j(K_j(\psi,\delta)))$ of the inner map in coordinates, and not by the position of the value $Q_j(\psi,\delta)$ of the return map in the plane to the left or right of the vertical axis. Our choice of disjoint sets  circumvents a discussion how the latter are related to the more accessible angles $\Phi_j(\psi,\delta)$.

\begin{proposition}
Assume the hypotheses of Proposition 7.3 are satisfied and let $\Delta_1=\Delta_1(\Delta_2)$. Consider the disjoint sets
$$
M_0= \{(\psi,\delta)\in[-\alpha_j,\alpha_j]\times[\Delta_1,\Delta_2]:\psi_j-\pi<\Phi_j(\psi,\delta)<\psi_j\}
$$
and
$$
M_1 = \{(\psi,\delta)\in[-\alpha_j,\alpha_j]\times[\Delta_1,\Delta_2]:\psi_j<\Phi_j(\psi,\delta)<\psi_j+\pi\}.
$$
For every curve $g:[a,b]\to[-\alpha_j,\alpha_j]\times[\Delta_1,\Delta_2]$ with $g_2(a)=\Delta_1$ and $g_2(b)=\Delta_2$ there exist $a_0<b_0<a_1<b_1$ in $[a,b]$ such that 
\begin{eqnarray*}
\mbox{on} &(a_0,b_0), & 
g(t)\in M_0\quad\mbox{and}\quad Q_{j,2}(g(t))\in(\Delta_1,\Delta_2),\\
& \mbox{with} & Q_{j,2}(g(a_0))=\Delta_1\quad\mbox{and}\quad Q_{j,2}(g(b_0))=\Delta_2,\\
\mbox{while on} &(a_1,b_1), & g(t)\in M_1\quad\mbox{and}\quad Q_{j,2}(g(t))\in(\Delta_1,\Delta_2),\\
& \mbox{with} & Q_{j,2}(g(a_1))=\Delta_2\quad\mbox{and}\quad Q_{j,2}(g(b_1))=\Delta_1.
\end{eqnarray*}
\end{proposition}

See Figure 6.

\medskip

\begin{figure}
	\includegraphics[page=1,scale=0.7]{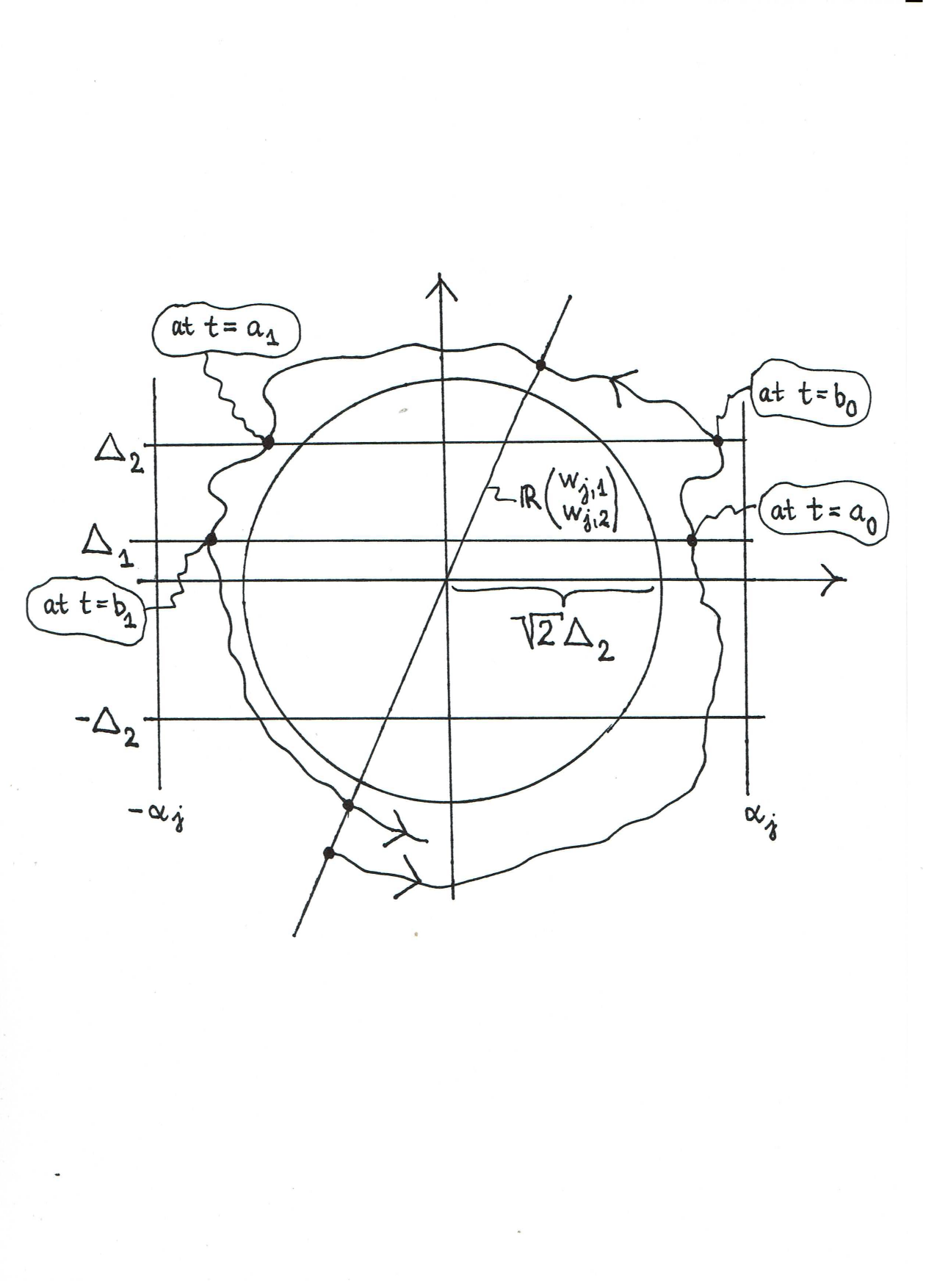}
	\caption{The values $Q_j(g(t))$ for $a\le t\le b$}
\end{figure}

{\bf Proof.} Proposition 7.2 yields $a_0'<b_0'\le a_1'<b_1'$ in $[a,b]$ such that
\begin{eqnarray*}
	\mbox{on} &(a_0',b_0'), & 
\Phi_j(g(t))\in(\psi_j-\pi,\psi_j),\\
& \mbox{with} & \Phi_j(g(a_0'))=\psi_j-\pi\quad\mbox{and}\quad \Phi_j(g(b_0'))=\psi_j,\\	\mbox{and on} &(a_1',b_1'), & \Phi_j(g(t))\in(\psi_j,\psi_j+\pi),\\
& \mbox{with} & \Phi_j(g(a_1'))=\psi_j\quad\mbox{and}\quad \Phi_j(g(b_1'))_2=\psi_j+\pi.
\end{eqnarray*}
From Proposition 7.3 (ii),
$$
Q_{j,2}(g(a_0'))<-\Delta_2,\quad Q_{j,2}(g(b_0'))>\Delta_2,\quad Q_{j,2}(g(a_1'))>\Delta_2,\quad Q_{j,2}(g(b_1'))<-\Delta_2.
$$
As in the proof of Proposition 7.2 one finds $a_0<b_0$ in $(a_0',b_0')$ and $a_1<b_1$ in $(a_1',b_1')$ with
\begin{eqnarray*}
Q_{j,2}(g(a_0))=\Delta_1 & \mbox{and} & Q_{j,2}(g(b_0))=\Delta_2,\quad\mbox{and}\quad Q_{j,2}(g(t))\in(\Delta_1,\Delta_2)\quad\mbox{on}\quad(a_0,b_0),\\
Q_{j,2}(g(a_1))=\Delta_2 & \mbox{and} & Q_{j,2}(g(b_1))_2=\Delta_1,\quad\mbox{and}\quad Q_{j,2}(g(t))\in(\Delta_1,\Delta_2)\quad\mbox{on}\quad(a_1,b_1).\\
\end{eqnarray*}
Observe that on $(a_0,b_0)\subset(a_0',b_0')$ we have $g(t)\in M_0$ while on $(a_1,b_1)\subset(a_1',b_1')$ we have $g(t)\in M_1$.
$\Box$

\section{Complicated dynamics}

For the results of this section we assume that the hypotheses of Proposition 7.4 are satisfied. It may be convenient to repeat all of these assumption here, beginning with the choice of a real number $\tilde{\eta}>0$ with
$$
\tilde{\eta}<\min\{\mu,-\sigma/2\}
$$ 
which satisfies the inequality (18).\\

\noindent
The numbers $\eta>0$, $j\in\mathbb{N}$, and $\delta_j\in(0,1)$ are chosen so that the relations (7-10) hold.\\

\noindent
For given $\beta\in(0,1/2]$, the reals $\alpha_j=\alpha_j(\beta)\in(0,\delta_j)$ and $\delta_{\beta,j}\in
(0,2\alpha_j/3]$ are chosen according to Proposition 6.2.\\ 

\noindent
$\Delta_2\in(0,\delta_{\beta,j})$ is chosen so that the inequality (19) holds, and $\Delta_1=\Delta_1(\Delta_2)$.\\

\medskip

Recall the disjoint sets $M_0$ and $M_1$ from Proposition 7.4.

\begin{proposition}
For every sequence $(s_n)_{n=0}^{\infty}$ in $\{0,1\}$ there are forward trajectories $(x_n)_{n=0}^{\infty}$ of  $Q_j$ in $[-\alpha_j,\alpha_j]\times[\Delta_1,\Delta_2]$ with $x_n\in M_{s_n}$ and $\quad\Delta_1\le Q_{j,2}(x_n)\le\Delta_2$ for all integers $\quad n\ge0$.
\end{proposition}

{\bf Proof.} 1. Let a sequence $(s_n)_{n=0}^{\infty}$ in $\{0,1\}$ be given. Choose a curve $g:[a,b]\to [-\alpha_j ,\alpha_j]\times[\Delta_1,\Delta_2]$ such that $g_2(t)\in(\Delta_1,\Delta_2)$ for $a<t<b$ and $g_2(a)=\Delta_1$, $g_2(b)=\Delta_2$, for example, $g(t)=(0,t)$ for $a=\Delta_1\le t\le\Delta_2=b$.

\medskip

For integers $n\ge0$ we construct recursively curves $g_n:[A_n,B_n]\to[-\alpha_j,\alpha_j]\times[\Delta_1,\Delta_2]$ with decreasing domains in $[a,b]$ as follows.

\medskip

1.1. In order to find $g_0$ we apply Proposition 7.4 to the curve $g$ and obtain $a_0<b_0<a_1<b_1$ in $[a,b]$ with the properties stated in Proposition 7.4. In case $s_0=0$ we define $g_0$ by $A_0=a_0$, $B_0=b_0$, $g_0(t)=g(t)$ for $A_0\le t\le B_0$. Notice that $g_0(t)\in M_{s_0}$ for all $t\in(A_0,B_0)$,  $Q_{j,2}(g_0(t))\in(\Delta_1,\Delta_2)$ on $(A_0,B_0)$,  $Q_{j,2}(g_0(A_0))=\Delta_1$, and $Q_{j,2}(g_0(B_0))=\Delta_2$. In case  $s_0=1$ we define $g_0$ by $A_0=a_1$, $B_0=b_1$, $g_0(t)=g(a_1+b_1-t)$ for $A_0\le t\le B_0$. Notice that also in this case $g_0(t)\in M_{s_0}$ for all $t\in(A_0,B_0)$,  $Q_{j,2}(g_0(t))\in(\Delta_1,\Delta_2)$ on $(A_0,B_0)$,  and $Q_{j,2}(g_0(A_0))=Q_{j,2}(g(a_1+b_1-a_1))=\Delta_1$, $Q_{j,2}(g_0(B_0))=Q_{j,2}(g(a_1+b_1-b_1))=\Delta_2$.

\medskip

1.2.  For an integer $n\ge0$ let a curve $g_n:[A_n,B_n]\to[-\alpha_j,\alpha_j]\times[\Delta_1,\Delta_2]$ be given with  $g_n(t)\in M_{s_n}$ for all $t\in(A_n,B_n)$ and $Q_{j,2}(g_n(t))\in(\Delta_1,\Delta_2)$ on $(A_n,B_n)$,  $Q_{j,2}(g_n(A_n))=\Delta_1$, $Q_{j,2}(g_n(B_n))=\Delta_2$. Proceeding as in Part 1.1, with the curve $[A_n,B_n]\ni t\mapsto Q_j(g_n(t))\in[-\alpha_j,\alpha_j]\times[\Delta_1,\Delta_2]$
in place of the former curve $g$, we obtain $A_{n+1}<B_{n+1}$ in $[A_n,B_n]$ and a curve $g_{n+1}:[A_{n+1},B_{n+1}]\to[-\alpha_j,\alpha_j]\times[\Delta_1,\Delta_2]$ with  $g_{n+1}(t)\in M_{s_{n+1}}$ for all $t\in(A_{n+1},B_{n+1})$ and $Q_{j,2}(g_{n+1}(t))\in(\Delta_1,\Delta_2)$ on $(A_{n+1},B_{n+1})$,  $Q_{j,2}(g_{n+1}(A_{n+1}))=\Delta_1$, $Q_{j,2}(g_{n+1}(B_{n+1}))=\Delta_2$.

\medskip

2. From $A_n\le A_{n+1}<B_{n+1}\le B_n$ for all integers $n\ge0$ we get $\cap_{n\ge0}[A_n,B_n]=[A,B]$ with $A=\lim_{n\to\infty}A_n\le\lim_{n\to\infty}B_n=B$. Choose $t\in[A,B]$ and define
$$
x_n=g_n(t)\in[-\alpha_j,\alpha_j]\times[\Delta_1,\Delta_2]
$$
for every integer $n\ge0$. This yields a forward trajectory of $Q_j$ with $\Delta_1\le Q_{j,2}(x_n)\le\Delta_2$ for all integers $n\ge0$.

Let an integer $n\ge0$ be given. Proof of $x_n\in M_{s_n}$. We have
$x_n=g_n(t)$ with $A_n\le A\le t\le B\le B_n$, and for all $v\in(A_n,B_n)$, $g_n(v)\in M_{s_n}$. By continuity, $x_n\in\, cl\,M_{s_n}$. In case $s_n=0$ this yields
$$
\psi_j-\pi\le\Phi_j(x_n)\le\psi.
$$
Assume $\Phi_j(x_n)\in\{\psi-\pi,\psi\}$, Then Proposition 7.3 (ii) gives $|Q_{j,2}(x_n)|>\Delta_2$, which contradicts the previous estimate of $Q_{j,2}(x_n)$. Consequently,
$$
\psi_j-\pi<\Phi_j(x_n)<\psi,\quad\mbox{or,}\quad x_n\in M_0=M_{s_n}.
$$
In case $s_n=1$ the proof is analogous. $\Box$  

\medskip

The final result extends Proposition 8.1 to  entire trajectories.
  
\begin{theorem}
For every sequence $(s_n)_{n=-\infty}^{\infty}$ in $\{0,1\}$ there exist entire trajectories $(y_n)_{n=-\infty}^{\infty}$ of $Q_j$ with
$y_n\in M_{s_n}$ for all integers $\,n$.
\end{theorem}

{\bf Proof.} 1. Let $(s_n)_{n=-\infty}^{\infty}$ in $\{0,1\}$ be given. Proposition 8.1 guarantees that for every integer $k$ there is a forward trajectory $(y_{k,n})_{n=0}^{\infty}$ of $Q_j$ in $[-\alpha_j,\alpha_j]\times[\Delta_1,\Delta_2]$ so that for each integer $n\ge0$,
$$
y_{k,n}\in M_{s_{n-k}}\quad\mbox{and}\quad\Delta_1\le Q_{j,2}(y_{k,n})\le\Delta_2.
$$ 
For integers $k,n$ with $k\ge-n$ we define
$$
z_{k,n}=y_{k,n+k},
$$
so that 
\begin{eqnarray*}
z_{k,n}=y_{k,n+k} & \in &  M_{s_{n+k-k}}=M_{s_n},\\ 
z_{k,n+1}=y_{k,n+1+k}=Q_j(y_{k,n+k}) & = & Q_j(z_{k,n}),\\
Q_{j,2}(z_{k,n})=Q_{j,2}(y_{k,n+k}) & \in & [\Delta_1,\Delta_2].
\end{eqnarray*}

1.1. Choice of subsequences for integers $N\ge0$.

\medskip

1.1.1. The case $N=0$: For every integer $k\ge0$, $z_{k,0}\in M_{s_0}$. The sequence $(z_{k,0})_{k=0}^{\infty}$ in the compact set 
$[-\alpha_j,\alpha_j]\times[\Delta_1,\Delta_2]$  has a convergent subsequence  $(z_{\kappa_0(k),0})_{k=0}^{\infty}$ given by a strictly increasing map $\kappa_0:\mathbb{N}_0\to\mathbb{N}_0$. Let  $y_0=\lim_{k\to\infty}z_{\kappa_0(k),0}\in\,cl\,M_{s_0}\subset[-\alpha_j,\alpha_j]\times[\Delta_1,\Delta_2]$. 

\medskip

1.1.2. The case $N=1$: We choose a convergent subsequence
of $(z_{\kappa_0(k),-1})_{k=0}^{\infty}$ given by a strictly increasing map $\mu_{-1}:\mathbb{N}_0\to\mathbb{N}_0$, and upon this a convergent subsequence of $(z_{\kappa_0\circ\mu_{-1}(k),1})_{k=0}^{\infty}$ given by
a strictly increasing map
$\mu_1:\mathbb{N}_0\to\mathbb{N}_0$, and define $\kappa_1=\mu_{-1}\circ\mu_{1}$. For $k\to\infty$,
$$
z_{\kappa_0\circ\kappa_1(k),-0},\quad z_{\kappa_0\circ\kappa_1(k),-1}\quad\mbox{and}\quad z_{\kappa_0\circ\kappa_1(k),1}
$$
converge to elements $y_0\in\,cl\,M_{s_0}$, $y_{-1}\in\,cl\,M_{s_{-1}}$, and $y_1\in\,cl\,M_{s_1}$, respectively.

1.1.3. The general case $N\in\mathbb{N}$: Consecutively choosing further convergent subsequences analogously to Part 1.1.2 we obtain strictly increasing maps $\kappa_n:\mathbb{N}_0\to\mathbb{N}_0$, $n=0,\ldots,N$,
so that for each $n\in\{-N,\ldots,N\}$ the sequence
$$
(z_{\kappa_0\circ\ldots\circ\kappa_N(k),n})_{k=0}^{\infty}
$$
converges for $k\to\infty$ to some $y_n\in\,cl\,M_{s_n}$.

\medskip

2. The {\it diagonal sequence} $K:\mathbb{N}_0\to\mathbb{N}_0$ defined by $K(N)=(\kappa_0\circ\ldots\circ\kappa_N)(N)$ is strictly increasing since for every $N\in\mathbb{N}_0$ we have
$$
K(N+1)=(\kappa_0\circ\ldots\circ\kappa_N)(\kappa_{N+1}(N+1))>
(\kappa_0\circ\ldots\circ\kappa_N)(\kappa_{N+1}(N))\ge
(\kappa_0\circ\ldots\circ\kappa_N)(N)=K(N)
$$
due to strict monotonicity of all maps involved.

\medskip

3. Let an integer $n$ be given and set $N=|n|$. Proof that 
$$
(z_{K(k),n})_{k=N+1}^{\infty}\quad\mbox{is a subsequence of}\quad
(z_{(\kappa_0\circ\ldots\circ\kappa_N)(k),n})_{k=N+1}^{\infty}.
$$
Consider the map $\lambda:\{k\in\mathbb{N}_0:k>N\}\to\mathbb{N}_0$ given by $
\lambda(k)=(\kappa_{N+1}\circ\ldots\circ\kappa_k)(k).$
For every integer $k\ge N+1$,  
$$
K(k)=(\kappa_0\circ\ldots\circ\kappa_N)(\lambda(k)),
$$
and $\lambda$ is strictly increasing because analogously to Part 2 we have
\begin{eqnarray*}
\lambda(k+1) & = & (\kappa_{N+1}\circ\ldots\circ\kappa_{k+1})(k+1)=(\kappa_{N+1}\circ\ldots\circ\kappa_k)(\kappa_{k+1})(k+1))\\
& \ge &(\kappa_{N+1}\circ\ldots\circ\kappa_k)(k+1)> (\kappa_{N+1}\circ\ldots\circ\kappa_k)(k)
=\lambda(k)
\end{eqnarray*}
for every integer $k\ge N+1$.

Being a subsequence of $(z_{(\kappa_1\circ\ldots\circ\kappa_{|n|})(k),n})_{k=N+1}^{\infty}$ the sequence $(z_{K(k),n})_{k=N+1}^{\infty}$ converges for $k\to\infty$ to $y_n\in\,cl\,M_{s_n}$.

\medskip

4. We show that $(y_n)_{n=-\infty}^{\infty}$ is an entire trajectory of $Q_j$. Let an integer $n$ be given and set $N=|n|$. From Part 3 in combination with Part 1.1.3 we get that
$(z_{K(k),n})_{k=N+1}^{\infty}$ converges to $y_n\in[-\alpha_j,\alpha_j]\times[\Delta_1,\Delta_2]$ and that
$(z_{K(k),n+1})_{k=N+2}^{\infty}$ converges to $y_{n+1}$.
Recall
$z_{k,n+1}=Q_j(z_{k,n})$ for all integers $k,n$ with $k\ge-n$. For integers $k>N=|n|$ we have $K(k)\ge k\ge-n$,
and the preceding statement yields
$$
z_{K(k),n+1}=Q_j(z_{K(k),n}).
$$
It follows that
$$
y_{n+1}=\lim_{N+2\le k\to\infty}z_{K(k),n+1}=\lim_{N+1\le k\to\infty}Q_j(z_{K(k),n})=Q_j(y_n).
$$

\medskip

5. Proof of $y_n\in M_{s_n}$ for all integers $n$. Let an integer $n$ be given. We have $y_n\in\,cl\,M_{s_n}$. In case $s_n=0$ this yields $\psi_j-\pi\le\Phi_j(y_n)\le\psi_j$. Therefore the assumption $y_n\notin M_{s_n}$ results in  $\Phi_j(y_n)\in\{\psi_j-\pi,\psi_j\}$, which according to Proposition 7.3 (ii) means $|Q_{j,2}(y_n)|>\Delta_2$, in contradiction to $\Delta_1\le Q_{j,2}(y_n)\le\Delta_2$. The proof in case $s_n=1$ is analogous. $\Box$

\section{Appendix: From the vectorfield $V$ to the flow $F$} 

Consider Shilnikov's scenario according to Section 1, with a  continuously differentiable vectorfield $V:\mathbb{R}^3\supset dom_V\to\mathbb{R}^3$, $dom_V$ open and $V(0)=0$, so that there is  a homoclinic solution $h_V$ of Eq. (1), and the eigenvalues $u>0$ and $\sigma\pm i\mu$, $\sigma<0<\mu$, of $DV(0)$ satisfy

\medskip

(H) $\quad\quad0<\sigma+u.$

\medskip

For simplicity assume in addition $dom_V=\mathbb{R}^3$ and that $V$ is bounded. This is not a severe restriction since we are only interested in solutions close to the compact homoclinic loop $cl\,h_V(\mathbb{R})=h_V(\mathbb{R})\cup\{0\}$, and one can achieve the desired properties by a modification of the vectorfield outside a neighbourhood of $cl\,h_V(\mathbb{R})$. - Then the solutions of Eq.(1) constitute a continuously differentiable flow $F_V:\mathbb{R}\times\mathbb{R}^3\to\mathbb{R}^3$.  In the sequel we describe how to transform $F_V$ into a flow $F$ with the properties (F1)-(F5) stated in Section 2.

\medskip

Choose eigenvectors $w\in\mathbb{R}^3$ for the eigenvalue $u>0$ of $DV(0)$, and
$z\in\mathbb{C}^3$ for the eigenvalue $\sigma+i\,\mu$ of the complexification of $DV(0)$. Set $w_1=\Re z\in\mathbb{R}^3$ and $w_2=\Im z\in\mathbb{R}^3$. Then $E_s=\mathbb{R}w_1\oplus\mathbb{R}w_2$ and $E_u=\mathbb{R}w$.
are invariant under $DV(0)$. For the isomorphism $B:\mathbb{R}^3\to\mathbb{R}^3$ given by  $Be_1=w_1,Be_2=w_2,Be_3=w$ we have that the linear map $B^{-1} DV(0) B$ is multiplication with the matrix
$$
A=\left(\begin{array}{lcr}\sigma & \mu & 0\\-\mu & \sigma & 0\\0 & 0 & u\end{array}\right),
$$
which leaves $L=B^{-1}E_s$ and $U=B^{-1}E_u$ invariant.
\medskip
 
Recall the local stable and unstable manifolds $W_s$ and $W_u$ of Eq. (1) at the origin. Given any $\lambda>0$ we may assume that with some $r=r(\lambda)>0$ they have the form 
$$
W_s=\{x+w_s(x):x\in E_s,|x|<r\}\quad\mbox{and}\quad W_u=\{y+w_u(y):y\in E_u,|y|<r\}
$$ 
for continuously differentiable maps 
$$
w_s:\{x\in E_s:|x|<r\}\to E_u\quad\mbox{and}\quad w_u:\{y\in E_u:|y|<r\}\to E_s
$$
which satisfy $w_s(0)=0,Dw_s(0)=0,w_u(0)=0,Dw_u(0)=0$, and 
$$
|Dw_s(x)|<\lambda\quad{for}\quad x\in E_s\quad\mbox{with}\quad|x|<r,\quad|Dw_u(y)|<\lambda\quad{for}\quad y\in E_u\quad\mbox{with}\quad|y|<r.
$$
For sufficiently small neigbourhoods $N$ of the origin in $\mathbb{R}^3$ the submanifolds $W_s$ and $W_u$ are invariant under $F_V$  in $N$, and $F_V(t,x)\in N$ for all $t\ge0$ implies $x\in W_s$ while $F_V(t,x)\in N$ for all $t\le0$ implies $x\in W_u$. 

\medskip

There exist $\hat{r}\in(0,r)$ and a continuously differentiable extension $\hat{w}_s:E_s\to E_u$ of the restriction of $w_s$ to $\{x\in E_s:|x|\le\hat{r}\}$  with $|D\hat{w}_s(x)|<\lambda$ on $E_s$, and analogously there is a continuously differentiable extension  $\hat{w}_ u:E_u\to E_s$ of the restriction of $w_u$ to $\{y\in E_u:|y|\le\hat{r}\}$ with $|D\hat{w}_u(y)|<\lambda$ on $E_u$. For $\lambda>0$ sufficiently small the continuously differentiable map 
$$
S:\mathbb{R}^3\to\mathbb{R}^3,\quad S(x+y)=x+y-\hat{w}_s(x)-\hat{w}_u(y)\quad\mbox{for}\quad x\in E_s,y\in E_u,
$$
with $S(0)=0$ and $DS(0)=id_{\mathbb{R}^3}$ satisfies
$$
|DS(z)-id_{\mathbb{R}^3}|\le\frac{1}{2}\quad\mbox{for all}\quad z\in\mathbb{R}^3.
$$
It follows that all derivatives of $S$ are isomorphisms. Using integration along line segments and the previous estimate one shows that $S$ is one-to-one.  In order to see that $S$ is onto notice that for every $\zeta\in\mathbb{R}^3$ the map $z\mapsto \zeta-(S(z)-z)$ is a strict contraction, whose fixed point is a preimage of $\zeta$ under $S$. As all  derivatives of $S$ are isomorphisms applications of the Inverse Mapping Theorem yield that $S^{-1}$ is continuously differentiable. Altogether, $S$ is a diffeomorphism which maps $\hat{W}_s=\{x+\hat{w}_s(x):x\in E_s\}$ onto $E_s$ and $\hat{W}_u=\{y+\hat{w}_u(y):y\in E_u\}$ onto $E_u$ .

\medskip 

Define the continuously differentiable flow $F:\mathbb{R}\times\mathbb{R}^3\to\mathbb{R}^3$ by
$$
F(t,x)=B^{-1}S(F_V(t,S^{-1}(Bx))).
$$
For all $t\in\mathbb{R}$, $F(t,0)=0$, which is property (F1). 

\medskip

On (F2): Necessarily, $h_V'(t)\neq0$ for all $t\in\mathbb{R}$. From the properties $h_V(t)\neq0\neq h'_V(t)$ for all $t\in\mathbb{R}$ and $\lim_{|t|\to\infty}h_V(t)=0$  in combination with the fact that $B^{-1}\circ S$ is a diffeomorphism with fixed point $0$ we infer that the flowline $h=B^{-1}\circ S\circ h_V$ of $F$ satisfies $h(t)\neq0\neq h'(t)$ for all $t\in\mathbb{R}$ and $\lim_{|t|\to\infty}h(t)=0$.

\medskip

On (F3): For $t\in\mathbb{R}$ let $T(t)=D_2F(t,0)$. Then
$$
T(t)=B^{-1} DS(0)D_2F_V(t,S^{-1}(B0))DS^{-1}(0)B=B^{-1}D_2F_V(t,0)B=B^{-1}e^{t\,DV(0)}B=e^{t\,B^{-1}DV(0)B}.
$$
As the linear map $B^{-1}DV(0)B$ is multiplication by the matrix $A$ we get $T(t)x=e^{t\,B^{-1}DV(0)B} x=e^{t\,A}\cdot x$ for all $x\in\mathbb{R}^3$ and $t\in\mathbb{R}$. It follows that $T(t)L\subset L$ and $T(t)U\subset U$  for every $t\in\mathbb{R}$.
The representation of $y=T(t)x=e^{t\,A}\cdot x$ described in property (F3) is confirmed by solving the initial value problems
$$
x'=Ax,\quad x(0)=e_{\nu},\quad\nu=1,2,3,
$$
whose solutions coincide at $t\in\mathbb{R}$ with the columns of the matrix $e^{t\,A}$.

\medskip

On (F4). We show that for a neighbourhood $N$ of the origin in $\mathbb{R}^3$ so small that $W_s$ is invariant under $F_V$ in $N$ and
$$
|x|<\hat{r}\quad\mbox{for all}\quad x\in E_s\quad\mbox{and}\quad y\in E_u\quad\mbox{with}\quad x+y\in N
$$ 
the space L is invariant under $F$  in $\tilde{N}=B^{-1}S(N)$.

\medskip

Proof. Let $z\in L\cap\tilde{N}$ and $t\in\mathbb{R}$ be given with $F(\tau,z)\in\tilde{N}$ for all $\tau\in t\cdot[0,1]$. As $B^{-1}\circ S$ maps $\hat{W}_s$ onto $L$ we get $x=S^{-1}(Bz)\in\hat{W}_s\cap N$, and for every $\tau\in t\cdot[0,1]$,
$$
F_V(\tau,x)=S^{-1}(BF(\tau,B^{-1}S(x)))=S^{-1}(BF(\tau,z))\in S^{-1}(B\tilde{N})=N.
$$
By the choice of $N$, $|\xi|<\hat{r}$ for $\xi\in E_s$ with $x=\xi+\hat{w}_s(\xi)$. Hence $\hat{w}_s(\xi)=w_s(\xi)$, and thereby $x\in W_s\cap N$. The invariance of $W_s$ under $F_V$ in $N$ yields $F_V(t,x)\in W_s\cap N$. As above, $|\xi_t|<\hat{r}$ for $\xi_t\in E_s$ with $F_V(t,x)=\xi_t+w_s(\xi_t)$. It follows that $w_s(\xi_t)=\hat{w}_s(\xi_t)$, and thereby
$$
F(t,z)=B^{-1}S(F_V(t,S^{-1}(Bz)))=B^{-1}S(F_V(t,x))=B^{-1}S(\xi_t+\hat{w}_s(\xi_t))\in L,
$$
which shows the desired invariance of $L$ under $F$ in $\tilde{N}$.

\medskip

It follows easily that for some $r_L>0$ the space $L$ is invariant  under $F$ in $\{z\in\mathbb{R}^3:|z|<r_L\}$. Analogously $U$ is invariant under $F$ in $\{z\in\mathbb{R}^3:|z|<r_U\}$, for some $r_U>0$. Let $r_F=\min\{r_L,r_U\}$. Then both $L$ and $U$ are invariant under $F$ in $\{z\in\mathbb{R}^3:|z|<r_F\}$, which is property (F4).

\medskip

On (F5). From $\lim_{t\to\infty}h_V(t)=0$, $h_V(t)\in W_s\cap\{x+y:x\in E_s,y\in E_u,|x|<\hat{r}\}$ for $t>0$ sufficiently large. By $\hat{w}_s(x)=w_s(x)$ for $x\in E_s$ with $|x|<\hat{r}$, 
$$
B^{-1}S(\{x+w_s(x):x\in E_s,|x|<\hat{r}\})\subset B^{-1}S(\hat{W}_s)=L.
$$ 
Consequently, there is $t_L\in\mathbb{R}$ with $h(t)=B^{-1}S(h_V(t))\in L$ for all $t\ge t_L$.
Analogously, $h(t)\in U$ for all $t\le t_U$, with some $t_u<t_L$. By continuity and $h(t)\neq0$ everywhere,
either $h(t)\in (0,\infty)e_3$ for all $t\le t_U$, or $h(t)\in(-\infty,0)e_3$ for all $t\le t_U$.

\end{document}